\documentclass{amsart}

\usepackage{amsmath,amssymb,url,color,mathtools}
\usepackage[margin=4cm]{geometry}
\usepackage{tikz}

\newtheorem{lemma}{Lemma}[section]
\newtheorem{prop}[lemma]{Proposition}
\newtheorem{cor}[lemma]{Corollary}
\newtheorem{thm}[lemma]{Theorem}
\newtheorem{example}[lemma]{Example}
\newtheorem{thm?}[lemma]{Theorem?}

\newtheorem{ques}[lemma]{Question}

\newtheorem{remark}[lemma]{Remark}

\begin{document}
\title[The Group-Theoretic Prime Ax-Katz Theorem]%
{Functional degrees and arithmetic applications II:
The Group-Theoretic Prime Ax-Katz Theorem}

\author{Pete L.\ Clark}
\author{Uwe Schauz}

\renewcommand\atop[2]{\genfrac{}{}{0pt}{}{#1}{#2}}
\newcommand{\Mod}[1]{\ (\mathrm{mod}\ #1)}
\newcommand{\etalchar}[1]{$^{#1}$}
\newcommand{\F}{\mathbb{F}}
\newcommand{\et}{\textrm{\'et}}
\newcommand{\ra}{\ensuremath{\rightarrow}}
\newcommand{\lra}{\ensuremath{\longrightarrow}}
\newcommand{\FF}{\F}
\newcommand{\ff}{\mathfrak{f}}
\newcommand{\Z}{\mathbb{Z}}
\newcommand{\N}{\mathbb{N}}
\newcommand{\NN}{\widetilde{\N}}
\newcommand{\mm}{\underline{m}}
\newcommand{\nn}{\underline{n}}
\newcommand{\ch}{}
\newcommand{\R}{\mathbb{R}}
\renewcommand{\P}{\mathbb{P}}
\newcommand{\PP}{\mathbf{P}}
\newcommand{\pp}{\mathfrak{p}}
\newcommand{\C}{\mathbb{C}}
\newcommand{\Q}{\mathbb{Q}}
\newcommand{\ab}{\operatorname{ab}}
\newcommand{\Aut}{\operatorname{Aut}}
\newcommand{\gk}{\mathfrak{g}_K}
\newcommand{\gq}{\mathfrak{g}_{\Q}}
\newcommand{\OQ}{\overline{\Q}}
\newcommand{\Out}{\operatorname{Out}}
\newcommand{\End}{\operatorname{End}}
\newcommand{\Gal}{\operatorname{Gal}}
\newcommand{\CT}{(\mathcal{C},\mathcal{T})}
\newcommand{\lcm}{\operatorname{lcm}}
\newcommand{\Div}{\operatorname{Div}}
\newcommand{\OO}{\mathcal{O}}
\newcommand{\rank}{\operatorname{rank}}
\newcommand{\tors}{\operatorname{tors}}
\newcommand{\IM}{\operatorname{IM}}
\newcommand{\CM}{\mathbf{CM}}
\newcommand{\HS}{\mathbf{HS}}
\newcommand{\Frac}{\operatorname{Frac}}
\newcommand{\Pic}{\operatorname{Pic}}
\newcommand{\coker}{\operatorname{coker}}
\newcommand{\Cl}{\operatorname{Cl}}
\newcommand{\loc}{\operatorname{loc}}
\newcommand{\GL}{\operatorname{GL}}
\newcommand{\PGL}{\operatorname{PGL}}
\newcommand{\PSL}{\operatorname{PSL}}
\newcommand{\Frob}{\operatorname{Frob}}
\newcommand{\Hom}{\operatorname{Hom}}
\newcommand{\Coker}{\operatorname{\coker}}
\newcommand{\Ker}{\ker}
\newcommand{\g}{\mathfrak{g}}
\newcommand{\sep}{\operatorname{sep}}
\newcommand{\new}{\operatorname{new}}
\newcommand{\Ok}{\mathcal{O}_K}
\newcommand{\ord}{\operatorname{ord}}
\newcommand{\Ohell}{\OO_{\ell^{\infty}}}
\newcommand{\cc}{\mathfrak{c}}
\newcommand{\ann}{\operatorname{ann}}
\renewcommand{\tt}{\mathfrak{t}}
\renewcommand{\cc}{\mathfrak{a}}
\renewcommand{\aa}{\mathbf{a}}
\newcommand\leg{\genfrac(){.4pt}{}}
\renewcommand{\gg}{\mathfrak{g}}
\renewcommand{\O}{\mathcal{O}}
\newcommand{\Spec}{\operatorname{Spec}}
\newcommand{\rr}{\mathfrak{r}}
\newcommand{\rad}{\operatorname{rad}}
\newcommand{\SL}{\operatorname{SL}}
\newcommand{\fdeg}{\operatorname{fdeg}}
\renewcommand{\rank}{\operatorname{rank}}
\newcommand{\Int}{\operatorname{Int}}
\newcommand{\zz}{\mathbf{z}}
\newcommand{\MaxSpec}{\operatorname{MaxSpec}}
\newcommand{\mmm}{\mathfrak{m}}

\begin{abstract}
We give a version of Ax-Katz's $p$-adic congruences and Moreno-Moreno's
$p$-weight refinement that holds over any finite commutative ring of prime
characteristic.   We deduce this from a purely group-theoretic result that gives a
lower bound on the $p$-adic divisibility of the number of simultaneous zeros of a
system of maps $f_j: A\to B_j$ from a fixed ``source'' finite commutative group $A$
of exponent $p$ to varying ``target'' finite commutative $p$-groups $B_j$.  Our proof
combines Wilson's proof of Ax-Katz over $\F_p$ with the functional calculus of
Aichinger-Moosbauer.
\end{abstract}

\maketitle


\section{Introduction}
\noindent This is the second in a sequence of papers in which we attempt a synthesis
and further development of work of Wilson \cite{Wilson06} and of Aichinger and
Moosbauer \cite{Aichinger-Moosbauer21}.  Whereas in the first paper
\cite{Clark-Schauz21} we applied arithmetic results of Weisman \cite{Weisman77}
and Wilson \cite{Wilson06} to answer a purely algebraic problem posed by
Aichinger-Moosbauer, in this paper the process is reversed: we use the algebraic
work of \cite{Clark-Schauz21} along with Aichinger-Moosbauer's functional calculus
to deduce arithmetic results.  In particular we give a purely group-theoretic result that
implies the theorem of Ax-Katz in the case of systems of polynomial equations over a
prime finite field $\F_p$ and the theorem of Moreno-Moreno on systems of polynomial
equations over a finite field $\F_q$.

\subsection{Notation and Terminology}
We denote by $\mathcal{P}$ the set of (positive) prime numbers,
write $\N$ for the set of non-negative integers, and put $\Z^+\!\coloneqq \N \setminus
\{0\}$. We endow the set \[ \NN \,\coloneqq\, \N \cup \{- \infty,\infty\} \] with the most
evident total ordering, in which $-\infty$ is the least element and $\infty$ is the
greatest element. The symbol $-\infty$ is also used as the degree of the zero
polynomial, which explains our restriction to nonzero polynomials or functions in some
theorems.
\\ \\
Throughout, $q = p^N$ denotes a positive integer power of a prime number $p$ and
$\F_q$ shall denote ``the'' (unique up to isomorphism) finite field of order $q$. For $n
\in \Z\setminus\{0\}$, we denote by $\ord_q(n)$ the largest power of $q$ that divides
$n$; we also put $\ord_q(0) = \infty$.
\\ \\
In this paper, rings are not necessarily commutative. We say that a ring $R$ is a
\textbf{domain} if for all $x,y \in R$, $xy = 0$ implies $x = 0$ or $y = 0$. A \textbf{rng}
is like a ring but not necessarily having a multiplicative identity. If $R,R_1, \dotsc,
R_r$ are sets, such that each of the sets $R_1, \dotsc, R_r$ contains a distinguished
element denoted $0$, and if $f_1:R^n\!\ra R_1$, \dots, $f_r:R^n\!\ra R_r$ are
functions (possibly given as polynomials), we also define
\[ Z(f_1,\dotsc,f_r)\,=\,Z_{R^n}(f_1,\dotsc,f_r)
 \,\coloneqq\, \bigl\{x \in R^n \mid f_1(x) = 0, \dotsc, f_r(x) = 0\bigr\}. \]

\subsection{Chevalley-Warning and Ax-Katz}
We begin by recalling the following results of Chevalley-Warning and Ax-Katz.

\begin{thm}\label{CWAK}
Let $p \in \mathcal{P}$ and $q:=p^N$\!. Let $f_1,\dotsc,f_r\in\F_q[t_1,\dotsc,t_n]$ be
nonzero polynomials. If $Z:=Z_{\F_q^n}(f_1,\dotsc,f_r)$ and $\sum_{j=1}^r \deg(f_j)<
n$\,, then
\begin{itemize}
\item[a)] $\,\ord_p(\# Z)\geq1$\hfill\emph{(Chevalley-Warning Theorem
    \cite{Chevalley35}, \cite{Warning35})},
\item[b)] $\,\ord_q(\# Z) \geq \Bigl{\lceil} \frac{n-\sum_{j=1}^r \deg(f_j)}{\max_{j=1}^r
    \deg(f_j)} \Bigr{\rceil}$\hfill\emph{(Ax-Katz Theorem \cite{Ax64}, \cite{Katz71})}.
\end{itemize}
\end{thm}
\noindent Theorem \ref{CWAK}b) in the case of one polynomial was
proved in 1964 by J.\ Ax \cite{Ax64}, while the general case was proved in 1971 by
N.M.\ Katz \cite{Katz71}.  Also in \cite{Ax64}, Ax gave a strikingly simple ten line
proof of Theorem \ref{CWAK}a).  There is certainly no known ten line proof of
Theorem \ref{CWAK}b): Ax's proof for one polynomial used methods of algebraic
number theory -- Jacobi sums and Stickelberger's congruence -- while Katz's proof of
the general case used some sophisticated arithmetic geometry -- zeta functions and
$p$-adic cohomology.  An Ax-style proof of Theorem \ref{CWAK}1.1b) was given by
D. Wan \cite{Wan89}, while Hou \cite{Hou05} gave a short deduction of Theorem
\ref{CWAK}1.1b) from the $r = 1$ case.  Also D.J.\ Katz \cite{Katz12} proved a result
in coding theory that implies Theorem \ref{CWAK}b).
\\ \\
What if we replace $\F_q$ by a finite ring $R$?  If $R$ is finite commutative and
\emph{principal} (i.e., every ideal of $R$ is principal), then for each prime number $p$
the largest power of $p$ dividing $\# Z_{R^n}(f_1,\dotsc,f_r)$ for all polynomials
$f_1,\dotsc,f_r \in R[t_1,\dotsc,t_n]$ of given positive degrees was determined: for $r
= 1$ by Marshall-Ramage \cite{Marshall-Ramage75} and in general by D.J.\ Katz
\cite{Katz09}.
\\ \indent
A finite commutative ring is Artinian, hence is a finite product of finite local Artinian
rings, each of which must have prime power order. In this way we immediately reduce
to the case of finite rings of prime power order.  Most such rings are however
\emph{not} principal, and there had been no known analogue of Chevalley-Warning --
let alone of Ax-Katz -- over any finite non-principal ring until the following recent
result.

\begin{thm}\emph{(Aichinger-Moosbauer \cite[Thm.\ 12.6]{Aichinger-Moosbauer21})}
\label{INTRO.0} Let $R$ be a finite rng of order a power of a prime number $p$, and
let $f_1,\dotsc,f_r\in R[t_1,\dotsc,t_n]$ be nonzero polynomials. If
$Z:=Z_{R^n}(f_1,\dotsc,f_r)$ and $\sum_{i=1}^r \deg(f_i) < n$\,, then
\[ \ord_p(\# Z) \,\geq\, 1. \]
\end{thm}

\noindent Our first main result gives a version of Ax-Katz for all finite rngs of
exponent $p$.

\begin{thm}[Ring-Theoretic Prime Ax-Katz Theorem]\label{AXKATZCOR}
Let $R$ be a finite rng with underlying additive group $(R,+)$ of prime exponent $p$,
so $(R,+) \cong \bigl((\Z/p\Z)^N\!,+\bigr)$ for some $N \in \Z^+$\!.  Let
$f_1,\dotsc,f_r\in R[t_1,\dotsc,t_n]$ be nonzero polynomials.  If
$Z:=Z_{R^n}(f_1,\dotsc,f_r)$, then
\begin{equation*}
\label{AXKATZCOREQ2}
\ord_p(\# Z) \,\geq\, \biggl{\lceil} \frac{N\bigl(n-\sum_{j=1}^r \deg(f_j)\bigr)}
{\max_{j=1}^r \deg(f_j)} \biggr{\rceil}.
\end{equation*}
\end{thm}

\begin{remark}
\label{CEILINGREMARK} If we take $R$ to be the finite field $\F_{p^N}$ of order
$p^N$\!, the conclusion of Theorem \ref{AXKATZCOR} is that
\begin{equation}
\label{AKREMARKEQ1}
\ord_p(\# Z) \,\geq\, \biggl{\lceil} \frac{N\bigl(n-\sum_{j=1}^r \deg(f_j)\bigr)}
{\max_{j=1}^r \deg(f_j)} \biggr{\rceil},
\end{equation}
while the Ax-Katz Theorem yields the $p$-adic congruence
\begin{equation}
\label{AKREMARKEQ2}
\ord_p(\# Z) \,\geq\, N\ord_{p^N}(\# Z)
 \,\geq\, N \biggl{\lceil} \frac{n-\sum_{j=1}^r \deg(f_j)}
 {\max_{j=1}^r \deg(f_j)} \biggr{\rceil}.
\end{equation}
The latter placement of the ceiling functions is more favorable, as the lower bound in
(\ref{AKREMARKEQ2}) is better than the lower bound in (\ref{AKREMARKEQ1}) if
$N>1$. This is why we speak of Theorem \ref{AXKATZCOR} as a generalization of
the ``Prime Ax-Katz Theorem'' and not of the Ax-Katz Theorem.
\end{remark}
\noindent Moreno-Moreno \cite{Moreno-Moreno95} used the Prime Ax-Katz Theorem
as input to give a different $p$-adic congruence for polynomial systems over any
finite field $\F_q$ that takes into account the $p$-weight degrees of the polynomials.
When $q > p$ the Moreno-Moreno $p$-adic congruences neither imply nor are
implied by the Ax-Katz $p$-adic congruences: cf. \cite[Thm.\ 0-1]{Moreno-Moreno95}.
In \S 4 we will give a $p$-weight version of Theorem \ref{AXKATZCOR} that
generalizes the Moreno-Moreno $p$-adic congruences from $\F_q$ to any finite
commutative ring of prime exponent.
\\ \\
Theorems \ref{INTRO.0} and \ref{AXKATZCOR} follow from deeper group-theoretic
results, as we now explain.

\subsection{The Aichinger-Moosbauer Functional Calculus}
In their recent work \cite{Aichinger-Moosbauer21}, Aichinger-Moosbauer developed a
fully fledged calculus of finite differences for functions $f: A \ra B$, where $A$ and
$B$ are commutative groups.   When $A$ and $B$ are $\R$-vector spaces, this
subject has a long pedigree, going back at least to work of Fr\'echet \cite{Frechet09}.
More recent works addressing the same topic include Leibman \cite{Leibman02} --
who works with not necessarily commutative groups -- and Laczkovich
\cite{Laczkovich04} -- who surveys and works to synthesize some of the prior
literature.  Neverthless, though the idea of such a calculus was not new,
Aichinger-Moosbauer's work is strikingly elegant, systematic and useful.
\\ \indent
We denote by $B^A$ the set of all functions $f: A \ra B$. It is a commutative group
under pointwise addition.  For each $a \in A$, we define a \textbf{difference operator}
$\Delta_a \in \End(B^A)$ by
\[\Delta_a f: x \longmapsto f(x+a)-f(x). \]
These endomorphisms all commute.  Following Aichinger-Moosbauer, we assign to
each
$f \in B^A$ a \textbf{functional degree} $\fdeg(f) \in \NN$ as follows: \\[4pt]
$\bullet$ We put $\fdeg(f) = -\infty$ if and only if $f =
0$.\footnote{Aichinger-Moosbauer in \cite{Aichinger-Moosbauer21} assign the
functional degree $0$ to the zero function. Here we follow the convention of
\cite{Clark-Schauz21}.  It certainly makes no
mathematical difference.}  \\[3pt]
$\bullet$ For $n \in \N$, we say that $\fdeg(f) \leq n$ if $\Delta_{a_1} \dotsm
\Delta_{a_{n+1}} f = 0$ for all $a_1,\dotsc,a_{n+1} \in A$.  If this holds for some
$n \in \N$, then $\fdeg(f)$ is the least $n$ for which it holds.  \\[3pt]
$\bullet$ If $\fdeg(f) \leq n$ holds for no $n \in \N$, then we put $\fdeg(f) = \infty$.
\\ \\
In other words, if we set $\sup(\emptyset):=-\infty$, then
\begin{equation}\label{eq.fdeg}
\fdeg(f)\,=\,\sup\bigl\{n\in\N\mid \exists a_1,\dotsc,a_{n} \in A,\,\Delta_{a_1} \dotsm
\Delta_{a_{n}} f\neq0\bigr\}\,.
\end{equation}
For commutative groups $A$ and $B$ and $d \in \N$, we put
\[ \mathcal{F}^d(A,B) \,\coloneqq\, \{f \in B^A \mid \fdeg(f) \leq d\}, \]
and we also put
\[ \mathcal{F}(A,B) \coloneqq \{f \in B^A \mid \fdeg(f) < \infty \}. \]
As introduced in \cite[\S 2]{Aichinger-Moosbauer21} and also discussed in \cite[\S
3]{Clark-Schauz21}, if $\Z[A]$ is the integral group ring of $A$, then the commutative
group $B^A$ has a canonical $\Z[A]$-module structure determined by the product
\[ [a]f:x\longmapsto[a]f(x) \coloneqq f(x+a)\]
of scalars of the form $[a]\in\Z[A]$ and vectors $f\in B^A$\!. In view of this, we may
equally well view $\Delta_a$ as the element $[a]-[0]$ of $\Z[A]$, since this element
acts on $B^A$ in the previously defined way.   We write $e(B)$ for the exponent
$\exp(B)$ if this number is finite and set $e(B):=0$ otherwise. This means
$e(B)\neq0$ if and only if there exists an $N \in \Z^+$ such that $Nb = 0$ for all $b \in
B$, and then $e(B)=\exp(B)$ is the least such $N$. With that definition, $B^A$ is
canonically a $\Z/e(B)\Z$-module, so we may also view $\Delta_a$ as living in the
group ring $(\Z/e(B)\Z)[A]$.
\\ \\
The functional degree gives a notion of ``polynomial function of degree $d$'' even
when there is no ring in sight. Moreover the notion of functional degree is partially
compatible with the degree of an actual polynomial function, in the following sense:

\begin{lemma}
\label{AM12.5} Let $R$ be a rng, let $f$ be a polynomial over $R$ in $n$ variables,
and let $E(f) \in R^{R^n}$ be the associated function. Then $\fdeg(E(f)) \leq \deg(f)$.
\end{lemma}
\begin{proof}
This is \cite[Lemma 12.5]{Aichinger-Moosbauer21}.
\end{proof}
\noindent Lemma \ref{AM12.5} shows that any discrepancy between
the functional degree and the degree of a polynomial map will only make
Chevalley-Warning\,/\,Ax-Katz type results stated in terms of the functional degree
\emph{stronger} than their classical analogues.
\\ \\
Here is the group-theoretic result of Aichinger-Moosbauer that underlies Theorem
\ref{INTRO.0}.

\begin{thm}[Group-Theoretic Chevalley-Warning Theorem]
\label{GTCW}
\textbf{} \\
Let $N,m,\alpha_1,\ldots,\alpha_m,n,\beta_1,\ldots,\beta_n,r \in \Z^+\!$, let $p \in
\mathcal{P}$, and let
\[ A \coloneqq \bigoplus_{i=1}^m \Z/p^{\alpha_i} \Z,\ B
\coloneqq \bigoplus_{i=1}^n \Z/p^{\beta_i}\Z \] be finite commutative $p$-groups. Let
$f_1,\dotsc,f_r: A^N\!\ra B$ be nonzero functions.  If $Z:=Z_{A^N}(f_1,\dotsc,f_r)$ and
\begin{equation*}
\label{GTCWEQ}
\biggl( \sum_{j=1}^r \fdeg(f_j) \biggr) \biggl( \sum_{i=1}^n (p^{\beta_i}\!-1) \biggr)
\,<\, \biggl( \sum_{i=1}^m p^{\alpha_i} -1 \biggr) N,
\end{equation*}
then
\[ \ord_p(\# Z)\,\geq\,1. \]
\end{thm}
\begin{proof}
This is \cite[Thm.\ 12.2]{Aichinger-Moosbauer21}.
\end{proof}
\noindent Applying Theorem \ref{GTCW} with $A = B = (R,+)$, the additive group of
a finite rng of order a power of $p$ and using Lemma \ref{AM12.5}, we deduce
Theorem \ref{INTRO.0}.
\\ \\
Here is the main result of this paper.

\begin{thm}
\label{AXKATZWILSON} Let $N,r,\beta_1,\dotsc,\beta_r  \in \Z^+$\!, let $p \in
\mathcal{P}$, and put $A := (\Z/p\Z)^N$\!. For each $1 \leq j \leq r$, let $f_j \in
(\Z/p^{\beta_j} \Z)^A$ be a nonzero function. If $Z:=Z_{A}(f_1,\dotsc,f_r)$, then
\[ \ord_p(\# Z) \geq \biggl\lceil \frac{N- \sum_{j=1}^r \frac{p^{\beta_j}-1}{p-1} \fdeg(f_j)}
 {\max_{j=1}^r \,p^{\beta_j-1} \fdeg(f_j)} \biggr\rceil. \]
\end{thm}

\begin{remark}\label{AXKATZWILSON-B}
The codomains of the maps $f_j$ in Theorem \ref{AXKATZWILSON} can easily be
generalized from cyclic $p$-groups $\Z/p^{\beta_j} \Z$ to arbitrary finite commutative
$p$-groups $B_j$. If, for each $1 \leq j \leq r$,
\[ B_j \,=\, \bigoplus_{k=1}^{K(j)} \Z/p^{\beta_{j,k} } \Z \quad\text{with}\ \,\beta_{j,1}
 \geq \dotsb \geq \beta_{j,K(j)} \geq 1\,, \]
then each of the given maps $f_j: A \ra B_j$ can be composed with the coordinate
projection $\pi_k: B_j \ra \Z/p^{\beta_{j,k} }\Z \eqqcolon B_{j,k}$\,, for $1\leq k\leq
K(j)$. This yields functions $f_{j,k} \coloneqq \pi_k \circ f_j$ with
$$\max_{1 \leq k \leq K(j)}\bigl(\fdeg(f_{j,k})\bigr)\,=\,\fdeg(f_j)\,.$$
Evidently, $f_j(x) = 0$ for all $j$ if and only if $f_{j,k}(x) = 0$ for all $j$ and $k$. So,
applying Theorem \ref{AXKATZWILSON} to the family of all maps $f_{j,k}: A \ra
B_{j,k}$ that are nonzero, we get
\[ \ord_p(\# Z(f_1,\dotsc,f_r)) \,\geq\,
 \biggl{\lceil} \frac{N - \sum_{j=1}^r \fdeg(f_j)\sum_{k=1}^{K(j)}\frac{p^{\beta_{j,k}}\!-1}{p-1}}
 {\max_{j=1}^r\,p^{\beta_{j,1}\!-1}\fdeg(f_j)}\biggr{\rceil}.\]
This result may be viewed as a generalization of Theorem \ref{AXKATZWILSON},
which we recover by taking each $B_j$ to be cyclic. In practice, however, this result
loses information from Theorem \ref{AXKATZWILSON} in that for each $j$ we use
only $\max_{1 \leq k \leq K(j)}\bigl(\fdeg( \pi_k \circ f_j)\bigr)$ instead of the individual
functional degrees of the maps $\pi_k \circ f_j$\,.
\end{remark}

\noindent We also have the following corollary, which generalizes Theorem
\ref{AXKATZCOR}:

\begin{cor}[Group-Theoretic Prime Ax-Katz Theorem]
\label{GTPAKT} Let $N,n,r \in \Z^+$\!, and put $A := (\Z/p\Z)^N$\!. Let $f_1,\dotsc,f_r
\in A^{A^n}$ be nonzero functions. If $Z:=Z_{A^n}(f_1,\dotsc,f_r)$, then
\[\ord_p(\# Z) \,\geq\,
\biggl{\lceil} \frac{N\bigl(n-\sum_{j=1}^r \fdeg(f_j)\bigr)}
{\max_{j=1}^r \fdeg(f_j)}\biggr{\rceil}.\]
\end{cor}

\begin{proof}
Let $\tilde{A} := A^n \cong (\Z/p\Z)^{n N}$.  For $1 \leq k \leq N$, let $\pi_k: A \ra
\Z/p\Z$ be the $k$th coordinate projection.  For $1 \leq j \leq r$ and $1 \leq k \leq N$,
put
 \[f_{j,k} \coloneqq \pi_k \circ f_j \in (\Z/p\Z)^{A^n} = (\Z/p\Z)^{\tilde{A}},
 \quad\text{with}\quad \fdeg(f_{j,k}) \,\leq\, \fdeg(f_j) \]
according to \cite[Lemma 3.8b)]{Clark-Schauz21}. For $x \in A^n$, we have $f_j(x) =
0$ for all $j$ if and only if $f_{j,i}(x) = 0$ for all $j$ and $i$. So, applying Theorem
\ref{AXKATZWILSON} to to the family of all maps $f_{j,k} \in (\Z/p^\Z)^{\tilde{A}}$ that
are nonzero, we get
 \[ \ord_p(\# Z)\,\geq\,\dotsb
 \,\geq\, \biggl{\lceil}\frac{Nn{-}{\textstyle\sum_{j=1}^r\sum_{k=1}^N \fdeg(f_j)}}
 {\max_{j=1}^r \fdeg(f_j)}\biggr{\rceil}
 \,=\,\biggl{\lceil} \frac{N\bigl(n-\sum_{j=1}^r \fdeg(f_j)\bigr)}{\max_{j=1}^r \fdeg(f_j)}
 \biggr{\rceil}. \qedhere \]
\end{proof}

\noindent If $R$ is a finite rng with underlying additive group $(R,+)$ finite of
exponent $p$, then applying Corollary \ref{GTPAKT} with $A = (R,+)$ and using
Lemma \ref{AM12.5}, we deduce Theorem \ref{AXKATZCOR}. Combining it instead
with a $p$-weight analogue of Lemma \ref{AM12.5} (Proposition
\ref{PWEIGHTPROP}), we will get our $p$-weight improvement of Theorem
\ref{AXKATZCOR} that recovers the Moreno-Moreno Theorem.

\begin{remark}
In an earlier version of our work, Corollary \ref{GTPAKT} was our main result, but
switching to Theorem \ref{AXKATZWILSON} made the proof \emph{easier}: cf.\
Remark \ref{AXKATZWILSON-B}.  The idea to this improvement arose from a draft
manuscript \cite{GGZ} that D.\ Grynkiewicz sent us in March of 2022.  These results are also contained
in the arxiv preprint \cite{Grynkiewicz22}. The statement
of our Theorem \ref{AXKATZWILSON} is directly inspired by \cite[Thm.\
1.3.22]{GGZ}, which is closely related to Theorem \ref{AXKATZWILSON} but involves
sums over residue systems modulo $p$ and reductions modulo powers of $p$ of
polynomials $f_1,\dotsc,f_r \in \Z[t_1,\dotsc,t_N]$ rather than arbitrary functions
between commutative $p$-groups.  Moreover, in a later draft of the same manuscript,
Grynkiewicz, Geroldinger and Zhong give a weighted version of their result.
\end{remark}

\subsection{Structure of the Paper}
\text{}\\[3pt]
$\bullet$ In \S 2 we give a canonical series representation for a map $f: A \ra
B$ between commutative groups of finite functional degree when $A$ is finitely
generated.  Moreover, for commutative domains of characteristic $0$, we explore the
connection between functions of finite functional degree and integer-valued
polynomials.
\\[3pt]
$\bullet$ In \S 3 we carry over a lemma of Wilson to our setting and then prove
Theorem \ref{AXKATZWILSON}.
\\[3pt]
$\bullet$ In \S 4 we discuss $p$-weights and prove a $p$-weight improvement of
Theorem \ref{AXKATZCOR}.
\\[3pt]
$\bullet$ In \S 5 we discuss work of the present authors \cite{Clark-Schauz23} and of Clark-Triantafillou \cite{Clark-Triantafillou23} that continues and complements the present work.


\subsection{Acknowledgments} Thanks to E. Aichinger for his interest in our present
work, which led to the communication of the results of
Geroldinger-Grynkiewicz-Zhong.  Thanks to D. Grynkiewicz for showing us two early
versions of \cite{GGZ}. Thanks to A.C. Cojocaru, N. Jones and N. Triantafillou for
stimulating conversations.

\section{The Fundamental Representation for $f \in B^{\Z^N}$}

\subsection{Preliminaries}

Let $N \in \Z^+$\!, and let $B$ be a commutative group.  In this section we give a
canonical series representation for functions $f \in \mathcal{F}(\Z^N\!,B)$ in terms of
\textbf{binomial polynomials}:
$$\textstyle\binom{t}{d} \coloneqq \frac{t(t-1)\dotsm(t-d+1)}{d!} \in\Q[t]
\quad\text{if}\quad d \in \Z^+\!.$$ Obviously, $\binom{x}{d}$ is an integer if $x \in \N$,
as it is the usual binomial coefficient. But, $\binom{x}{d}$ is always an integer, also for
negative $x\in\Z$: see e.g.\ \cite[p.\ 19]{Cahen-Chabert}. The binomial polynomials
$\binom{t}{d}\in\Q[t]$ are \textbf{integer-valued polynomials} as they give rise to
functions $\binom{x}{d}$ from $\Z$ to $\Z$. We also take $\binom{x}{0}: \Z \ra \Z$ to
be the constant function $1$. And, we define $\binom{x}{d}: \Z \ra \Z$ to be the zero
function for negative $d\in\Z$. We discuss this kind of functions in \S 2.3.
\\ \\
For $1 \leq i \leq n$, let $e_i$ be the $i$th standard basis vector of $\Z^N$\!.  We
write $\Delta_i$ for the difference operator $\Delta_{e_i}$ of $B^{\Z^N}$\!\!.

\begin{lemma}
\label{4.0} Let $\underline{B}$ be a subgroup of the commutative group $B$, and let
$f \in B^{\Z^N}$\!\!. Then the following properties are equivalent:
\begin{itemize}
\item[(i)] $f(\Z^N) \subseteq \underline{B}$.
\item[(ii)] $\Delta_i f(\Z^N) \subseteq \underline{B}$ for all $1 \leq i \leq N$, and
    $f(\underline{0}) \in \underline{B}$.
\end{itemize}
\end{lemma}
\begin{proof}
(i) $\Rightarrow$ (ii) is immediate.\smallskip  \\
(ii) $\Rightarrow$ (i): For any $\underline{x} \in \Z^N$ and any $1 \leq i \leq N$, we
have
\[ \Delta_i f(\underline{x}) \,=\, f(\underline{x}+e_i)-f(\underline{x})
\,\in\, \underline{B}, \]
which shows that $f(\underline{x}+e_i) \in \underline{B} \iff f(\underline{x}) \in
\underline{B}$.  Since $f(\underline{0}) \in B$, an immediate inductive argument now
shows that $f(\underline{x}) \in \underline{B}$ for all $\underline{x} \in \Z^N$.
\end{proof}
\noindent For $\nn \coloneqq (n_1,\dotsc,n_N) \in \N^N$\!, we put
\[ \Delta^{\nn} \,\coloneqq\, \Delta^{n_1}_1 \dotsm \Delta^{n_N}_N
 .\]
Because $e_1,\dotsc,e_N$ is a set of generators for $\Z^N$\!, the following
characterization of the functional degree follows from  or \cite[Lemmas
2.2]{Aichinger-Moosbauer21}, or from \cite[Lemma 3.11]{Clark-Schauz21}:

\begin{prop}\label{thm.fdeg}
Let $f \in B^{\Z^N}$\! and define $\sup(\emptyset):=-\infty$. Then
$$\fdeg(f)\,=\,\sup\bigl\{|\nn|\mid\nn\in\N^N\!,\,\Delta^{\nn}f\neq0\bigr\}\,.$$
\end{prop}

\noindent If we compare this expression with the following definition of the
\textbf{$j$-th partial functional degree} for functions $f$ in $B^{\Z^N}$\!\!, which is
given by
\begin{equation}
\label{eq.pfdeg}
\fdeg_j(f)\,:=\,\sup\bigl\{n\in\N\mid \Delta_j^{n}f\neq0\bigr\}\,,
\end{equation}
it is easy to see that, for each $1 \leq j \leq N$,
\begin{equation*}
\fdeg_j(f) \,\leq\, \fdeg(f) \,\leq\, \sum_{i=1}^N \fdeg_i(f).
\end{equation*}

All this can easily be generalized to domains that are a direct product of arbitrary
commutative groups $A_1,\dotsc,A_N$, as in \cite[\S 5]{Aichinger-Moosbauer21}.
Regarding this product as an internal direct product, we define the $j$-th partial
functional degree of a function $f \in B^{\bigoplus_{i=1}^N A_i}$ by
$$\fdeg_j(f)\,:=\,\sup\bigl\{n\in\N\mid \exists\,a_{1},\dotsc,a_{n}\in A_j,\,
 \Delta_{a_{n}}\dotsm\Delta_{a_{1}}f\neq0\bigr\}\,.$$
It follows from \cite[Lemmas 2.2]{Aichinger-Moosbauer21} or \cite[Lemma
3.11]{Clark-Schauz21} that
\begin{equation}
\label{eq.fdegn}
\fdeg(f)\,=\,\sup\bigl\{n\in\N\mid \exists\,a_{1},\dotsc,a_{n}\in A_1\cup\dotsb\cup A_n,\,
 \Delta_{a_{n}}\dotsm\Delta_{a_{1}}f\neq0\bigr\}\,
\end{equation}
From this we easily get \cite[Theorem 5.2]{Aichinger-Moosbauer21}, for which we
present a shortened proof:

\begin{thm}\label{thm.pfdeg}
Let $A_1,\dotsc,A_N,B$ be commutative groups, and let $f \in B^{\bigoplus_{i=1}^N
A_i}$. Then, for each $1 \leq j \leq N$,
\begin{equation*}
\fdeg_j(f) \,\leq\, \fdeg(f) \,\leq\, \sum_{i=1}^N \fdeg_i(f).
\end{equation*}
\end{thm}
\begin{proof}
We may assume $f\neq0$, as the inequality holds otherwise. Comparing
\eqref{eq.fdeg} and \eqref{eq.pfdeg}, we see that $\fdeg_j(f)\leq\fdeg(f)$. To prove
$\fdeg(f) \leq \sum_{i=1}^N \fdeg_i(f)=:n\geq0$, let $a_{1},\dotsc,a_{n+1}\in
A_1\cup\dotsb\cup A_n$. By \eqref{eq.fdegn}, it suffices to show that
$\Delta_{a_{n+1}}\dotsm\Delta_{a_{1}}f=0$. As $n+1>\sum_{i=1}^N \fdeg_i(f)$, there
exists a $1\leq j\leq N$ such that more than $n_j:=\fdeg_j(f)$ of the elements
$a_{1},\dotsc,a_{n+1}$ lie inside $A_j$. Without loss of generality, assume
$a_{1},\dotsc,a_{n_j+1}\in A_j$. Then $\Delta_{a_{n_j+1}}\dotsm\Delta_{a_{1}}f=0$,
by \eqref{eq.pfdeg}, and $\Delta_{a_{n+1}}\dotsm\Delta_{a_{1}}f=0$ follows.
\end{proof}

\noindent  For the convenience of the reader, we also restate \cite[Lemma
2.2]{Clark-Schauz21}.

\begin{lemma}
\label{LEMMA2.2} Let $A$ and $B$ be commutative groups.  Let $a \in A$, $n \in \N$
and let $\Delta_a^n$ be the $n$-fold product $\Delta_a \dotsm \Delta_a \in \End
B^A$\!. For all $f \in B^A$ and all $x \in A$, we have
\[ \Delta_a^n f(x) \,=\, \sum_{i=0}^n (-1)^i \tbinom{n}{i} f(x+(n-i)a)
\,=\, \sum_{j=0}^n (-1)^{n-j} \tbinom{n}{j} f(x+ja). \]
\end{lemma}
\noindent
We recall an old result for comparison and future use:

\begin{lemma}
\label{CATS} Let $R$ be a commutative domain, let $f \in R[t_1,\ldots,t_n]$, and let
$X_i$ be a nonempty subset of $R$ with $\# X_i > \deg_i(f)$, for each $1 \leq i \leq
n$.  If $f(x) = 0$ for all $x \in X:=\prod_{i=1}^n X_i$\,, then $f = 0$.
\end{lemma}
\begin{proof}
We can immediately reduce to the case in which $\# X_i = \deg_i(f) + 1$ for all $i$.
Then the case $R = \Z$ is \cite[Lemma 2.1]{Alon-Tarsi92}, and their proof works
verbatim over any commutative domain. More general results appear in \cite[\S
2]{Schauz08}; see also \cite[Thm. 12]{Clark14}.
\end{proof}

\noindent The following result is related to Lemma \ref{CATS} and also to
\cite[Thm.\ 2.5]{Schauz14}.

\begin{lemma}
\label{lem.0} Let $N \in \Z^+$\!, let $B$ be a commutative group and let $f \in
B^{\Z^N}$\!\!. For each $1 \leq i \leq N$, let $a_i\in\Z$, $d_i\in\N$ with
$d_i\geq\fdeg_i(f)$, and put  \[[a_i,a_i+d_i]:=\{a_i,a_i+1,\dotsc,a_i+d_i\}. \]
If $f(\underline{x}) = 0$ for all $\underline{x}\in\prod_{i=1}^N[a_i,a_i+d_i]$\,, then $f =
0$.
\end{lemma}

\begin{proof}
We proceed by induction on $N$.\smallskip  \\
\emph{Base Case:} Suppose that $N = 1$, i.e.\ $f \in B^{\Z}$\!, $\fdeg(f) \leq d_1$
and $f(a) = f(a+1) = \dotsb = f(a+d_1) = 0$.  Applying Lemma \ref{LEMMA2.2} with
$d_1+1$ in the place of that lemma's $n$, with $1$ in the place of that lemma's $a$,
and with $a-1$, resp.\ $a$, in the place of that lemma's $x$, we can deduce $f(a-1) =
0$, rsp.\ $f(a+d_1+1) = 0$. Repeating this argument
we get $\dotsb=f(a-2) = f(a-1) = 0$ and $0=f(a+d_1+1) = f(a+d_1+2) =\dotsb$, i.e.\ $f = 0$.
\smallskip\\
\emph{Induction Step:} Suppose that $N \geq 2$ and that the result holds for all $f \in
\mathcal{F}(\Z^{N-1}\!,B)$.  For $0 \leq j \leq d_N$, put
\[ g_j \coloneqq f(\,\cdot,\dotsc,\cdot,a_N +j): \Z^{N-1} \!\ra B. \]
Then we have $\fdeg_i g_j \leq d_i$ for all $1 \leq i \leq N-1$ and $g_j$ vanishes
identically on $\prod_{i=1}^{N-1} [a_i,a_i+d_i]$, so induction gives $g_j = 0$ for all $0
\leq j \leq d_N$.  It follows that for each fixed $(x_1,\dotsc,x_{N-1}) \in \Z^{N-1}$ the
function $f(x_1,\dotsc,x_{N-1},\cdot): \Z \ra B$ vanishes on $[a_N,a_N+d_N]$, and it
has functional degree at most $d_N$.  So, by the base case, these functions are
identically zero, which means that $f$ is identically zero.
\end{proof}


\begin{lemma}
\label{lem.0b} Let $B$ be a commutative group, and let $f:\N^N\!\longrightarrow B$
be a function. If $\Delta^{\nn}f(\underline{0}) = 0$ for all $\nn \in \N^N$\!, then $f$ is
the zero function.
\end{lemma}

\begin{proof}
Given a function $f:\N^N\!\ra B$ with $\Delta^{\nn}f(\underline{0}) = 0$ for all $\nn \in
\N^N$, we prove the formally stronger conclusion $\Delta^{\mm}f(\underline{x}) = 0$
for all $\underline{x},\mm\in \N^N$\!. We do this by induction on $|\underline{x}|
\coloneqq x_1 + \dotsb + x_N$.\smallskip  \\
\emph{Base Case:} If $|\underline{x}| = 0$ then $\underline{x} = \underline{0}$ and
$\Delta^{\mm}f(\underline{x}) = \Delta^{\mm}f(\underline{0}) = 0$ for all $\mm\in
\N^N$\!, by the assumption on $f$.\smallskip\\
\emph{Induction Step:} Let $\underline{0}\neq\underline{x}\in \N^N$ and assume that
the statement holds for all $\underline{z}\in \N^N$ with
$|\underline{z}|<|\underline{x}|$, and for all $\mm\in \N^N$\!. As
$\underline{x}\neq\underline{0}$, there is an index $i$ such that $x_i\geq 1$. Hence,
by the induction hypothesis, for each $\mm\in \N^N$\!,
\[\Delta^{\mm}f(\underline{x}-e_i)\,=\,0\quad\text{and}\quad
\Delta^{\mm+e_i}f(\underline{x}-e_i)\,=\,0,\] so that
\begin{align*}
\Delta^{\mm}f(\underline{x})
&\,=\,\Delta^{\mm}f(\underline{x})-\Delta^{\mm}f(\underline{x}-e_i)\\
&\,=\,\Delta_i(\Delta^{\mm}f)(\underline{x}-e_i)\\
&\,=\,\Delta^{\mm+e_i}f(\underline{x}-e_i)\\
&\,=\,0,
\end{align*}
completing the induction step and the proof.
\end{proof}

%

\subsection{The Fundamental Representation}
We can now prove the following result, on which much of the rest of this work is
based.

\begin{thm}
\label{4.2} Let $B$ be a commutative group, and let $f \in B^{\Z^N}$\!\!.
\begin{itemize}
\item[a)] There is a unique function $a_{\bullet}: \N^N\!\ra B$ such that
\begin{equation*}
\label{4.2EQ1}
f(\underline{x}) = \sum_{\nn  \in \N^N} \binom{x_1}{n_1} \dotsm
\binom{x_N}{n_N} a_{\nn} \quad\text{for all $\underline{x} \in \N^N$\!.}
\end{equation*}
The function values of $a_{\bullet}$ are given by the formula $a_{\nn} =
\Delta^{\nn}f(\underline{0})$. 
\smallskip
\item[b)] If $d:=\fdeg(f)< \infty$, then
\begin{equation*}
f(\underline{x}) \,= \sum_{\atop{\nn \in \N^N}{\!\!|\nn| \leq d}} \binom{x_1}{n_1}
\dotsm \binom{x_N}{n_N} \Delta^{\nn}f(\underline{0})\quad\text{for all
$\underline{x} \in \Z^N$\!.}
\end{equation*}
\end{itemize}
\end{thm}
\begin{proof}
a) To prove the uniqueness, assume there is an $a_{\bullet}$ with
$f(\underline{x}):=\sum \binom{x_1}{n_1} \dotsm \binom{x_N}{n_N}a_{\nn}$ for all
$\underline{x}\in\N^N$\!.  For each $\underline{x} = (x_1,\dotsc,x_N) \in \N^N$ we
have $\binom{x_1}{n_1} \dotsm \binom{x_N}{n_N} = 0$ unless $n_i \leq x_i$ for all $1
\leq i \leq N$, so for each fixed $\underline{x}$ we have a finite sum.  For all $n \in
\Z^+$\!, we have $\binom{x+1}{n} - \binom{x}{n} = \binom{x}{n-1}$. From this it follows
that, for all $\mm,\nn \in \N$,
\begin{equation*}
\Delta^{\mm}\left(\binom{x_1}{n_1} \dotsm \binom{x_N}{n_N}\right)(\underline{0})
\,=\, \prod_{i=1}^N \binom{0}{n_i-m_i} \,=\, \begin{cases} 1 & \text{if }\mm = \nn\,, \\ 0 &
\text{otherwise.} \end{cases}
\end{equation*}
With that, if we apply $\Delta^{\mm}$ to $f(\underline{x})=\sum \binom{x_1}{n_1}
\dotsm \binom{x_N}{n_N}a_{\nn}$ and evaluate at $\underline{0}$\,, we see that
$a_{\mm} = \Delta^{\mm}f(\underline{0})$. So, the function
$a_{\bullet}:\nn\mapsto\Delta^{\nn}f(\underline{0})$ is the only possible choice. \\
\indent To show that this choice indeed yields the function $f$, define $\hat f:\N^N\!\ra
B$ by
\[\hat f(\underline{x})\,:=\,\sum_{\nn  \in \N^N}
 \binom{x_1}{n_1}\dotsm\binom{x_N}{n_N}\Delta^{\nn}f(\underline{0}). \]
By what we have already proven about the uniqueness of coefficients, for each $\nn
\in \N^N$\!, the coefficient $\Delta^{\nn}f(\underline{0})$ must be equal to
$\Delta^{\nn}\hat f(\underline{0})$, i.e.\ $\Delta^{\nn}(f-\hat f)(\underline{0})=0$. With
that, Lemma \ref{lem.0b} yields $f-\hat f=0$, i.e.\ $f=\hat f$, as desired.
\smallskip  \\
b) We have $\Delta^{\nn} f(\underline{0}) = 0$ for all $\nn \in \N^N$ with $|\nn| > d$,
so
\begin{equation*}
f(\underline{x}) \,= \sum_{\atop{\nn \in \N^N}{\!\!|\nn| \leq d}} \binom{x_1}{n_1}
\dotsm \binom{x_N}{n_N} \Delta^{\nn}f(\underline{0})\quad\text{for all
$\underline{x} \in \N^N$\!.}
\end{equation*}
The right hand side of this equation, however, defines a function $P: \Z^N\!\ra B$.
And, $f-P$ has functional degree at most $d$ and vanishes on $\N^N$\!. So, by
Lemma \ref{lem.0}, $f = P$.
\end{proof}

\noindent A finite linear combination $\sum_{|\nn| \leq d} \binom{x_1}{n_1} \dotsm
\binom{x_N}{n_N} a_{\nn}$ of multivariate binomial polynomials $\binom{x_1}{n_1}
\dotsm \binom{x_N}{n_N}$ with coefficients $a_{\nn}$ in the group $B$, as in
Theorem \ref{4.2}b), was called a \textbf{polyfract} in \cite{Schauz14}. Due to the
uniqueness of the coefficients $a_{\nn}$ in Theorem \ref{4.2}a), we do not have to
distinguish between a polyfract $\sum_{|\nn| \leq d} \binom{x_1}{n_1} \dotsm
\binom{x_N}{n_N} a_{\nn}$\,, where the $x_i$ may be seen as symbolic variables,
and the corresponding polyfractal function
$$f:\Z^N\!\longrightarrow B\,,\ x\longmapsto\sum_{\atop{\nn \in \N^N}{\!\!|\nn| \leq d}} \binom{x_1}{n_1}
 \dotsm \binom{x_N}{n_N} a_{\nn}\,.$$
It is also clear that the functional degree of such a map is given by
\begin{equation}\label{eq.fdeg0}
\fdeg(f)\,=\,\sup\bigl\{|\nn| \leq d\mid\nn\in\N^N\!,\,a_{\nn}\neq0\bigr\}\,,
\end{equation}
where $\sup(\emptyset):=-\infty$. This follows from the fact that
$$\Delta_i^{m_i}\binom{x_i}{n_i} \,=\, \binom{x_i}{n_i-m_i}$$
and the observation that $x_i\mapsto\binom{x_i}{n_i-m_i}$ is not the zero map
whenever $m_i\leq n_i$.
%
%
So if we combine Theorem \ref{4.2}b) with formula \eqref{eq.fdeg0}, we obtain the
following corollary, which allows us to calculate the functional degree in a more
localized fashion -- at least if the functional degree is known to be finite.

\begin{cor}\label{4.2cor-d}
Let $B$ be a commutative group, and let $f \in B^{\Z^N}$\!\!. If $\fdeg(f)<\infty$ then
$$\fdeg(f)\,=\,\fdeg_{\underline{0}}(f)\,:=\,\sup\bigl\{|\nn|\mid\nn \in
    \N^N\!,\,\Delta^{\nn}f(\underline{0})\neq0\bigr\}\,.$$
\end{cor}

\noindent In this corollary, the point $\underline{0}$ can be replaced by any other
point $\underline{a} \in \Z^N\!$, since $\fdeg([\underline{a}]f)=\fdeg(f)$. If
$\fdeg(f)=\infty$, however, we may have $\fdeg_{\underline{0}}(f)<\infty$. This is
because $\fdeg_{\underline{0}}(f)$ depends only on the function values
$f(\underline{x})$ at points $\underline{x}\in\N^N$\!, whereas the requirement
$\fdeg(f)<\infty$ allows only one unique extension from $\N^N$ to $\Z^N$ -- the
extension given by the formula in Theorem \ref{4.2}b).\smallskip

\noindent We also see that in the case $B=\Q$, the series representation in Theorem
\ref{4.2}b) provides a polynomial $\hat f\in\Q[t_1,\dotsc,t_N]$ that describes $f$:

\begin{cor}\label{4.2cor}
If $f \in\Q^{\Z^N}$\! has finite functional degree, then there exists a polynomial $\hat
f\in\Q[t_1,\dotsc,t_N]$ with $\deg(\hat f)=\fdeg(f)$ and $\hat
f(\underline{x})=f(\underline{x})$ for all $\underline{x} \in\Z^N$\!.
\end{cor}


\begin{remark}\leavevmode
\begin{itemize}
\item[a)] For $B$ a finitely generated commutative group, the series representation
    in Theorem \ref{4.2}b) was explored in \cite[\S 2]{Schauz14}.\smallskip
\item[b)] The series expansion of Theorem \ref{4.2} is a discrete analogue of the
    Taylor series expansion of a smooth function $f: \R^N \!\ra \R$.  Theorem
    \ref{4.2}a) implies a \emph{uniqueness} property: for any two functions
    $a_{\bullet},b_{\bullet}: \N^N \!\ra B$ that each map all but finitely many elements
    of the domain to $0$, define associated functions
\[ f_{a_{\bullet}}: \Z^N \!\ra B, \ \underline{x} \mapsto \sum_{\nn\in \N^N}
\binom{x_1}{n_1}
\dotsm \binom{x_N}{n_N} a_{\nn} \]
and
\[ f_{b_{\bullet}}: \Z^N \!\ra B, \ \underline{x} \mapsto \sum_{\nn\in \N^N}
\binom{x_1}{n_1}
\dotsm \binom{x_N}{n_N} b_{\nn}. \] Then $f_{a_{\bullet}} = f_{b_{\bullet}}$ if and
only if $a_{\bullet} = b_{\bullet}$.  This is a discrete analogue of the fact that in a
power series expansion centered at $0$, the coefficients are determined by the
partial derivatives at $0$.\smallskip
\item[c)] Just as it is immediate to also consider Taylor series expansions centered
    at a nonzero point $a \in \R^N$\!, there are also representations of $f \in
    \mathcal{F}(\Z^N\!,B)$ based on the values $\Delta^{\nn} f(\underline{a})$ for any
    fixed $\underline{a} \in \Z^N$\!.
\end{itemize}
\end{remark}

\noindent Next we recall some notation and a result from \cite[\S
3.2]{Clark-Schauz21}.  If $\varepsilon: A\ra A$ and $\mu: B\ra B'$ are
homomorphisms of commutative groups, then we have group homomorphisms
\[ \varepsilon^*: B^A\!\longrightarrow B^{A'}\!,\, \
f\longmapsto \varepsilon^* f \coloneqq f \circ \varepsilon \]
and
\[ \mu_*: B^A\!\longrightarrow (B')^A,\, \ f\longmapsto \mu_* f \coloneqq \mu \circ f. \]
It is easy to see that $\varepsilon^*$ is injective if and only if $\varepsilon$ is surjective and that $\mu_*$ is surjective if and only if $\mu$ is surjective.
\\ \\
The following result is \cite[Lemma 3.9]{Clark-Schauz21}.

\begin{lemma}[Homomorphic Functoriality I]
\label{HOMFUNCI} Let $\varepsilon: A' \ra A$ and $\mu: B \ra B'$ be
homomorphisms of commutative groups, and let $f \in B^A$\!.  Then:
\begin{itemize}
\item[a)] $\fdeg \varepsilon^* f \leq \fdeg f$, with equality if $\varepsilon$ is
    surjective;
\item[b)] $\fdeg \mu_* f \leq \fdeg f$, with equality if $\mu$ is injective.
\end{itemize}
\end{lemma}

\noindent The following conceptually similar result is a consequence of Theorem
\ref{4.2}.

\begin{cor}[Homomorphic Functoriality II]
\label{HOMFUNCII} Let $B,B'$ be commutative groups, let $\mu: B \ra B'$ be a
homomorphism,  and let $f \in B^{\Z^N}$\!\!.
\begin{itemize}
\item[a)]\ \vspace{-.5em}
\begin{equation*}
\mu_*f(\underline{x}) \,= \sum_{\nn \in \N^N} \binom{x_1}{n_1}
\dotsm \binom{x_N}{n_N} \mu(\Delta^{\nn}f(\underline{0}))
\quad\text{for all $\underline{x} \in \N^N$\!.}
\end{equation*}
\item[b)] If $d:=\fdeg(\mu_*f)< \infty$, then
\begin{equation*}
\mu_*f(\underline{x}) \,= \sum_{\atop{\nn \in \N^N}{\!\!|\nn| \leq d}} \binom{x_1}{n_1}
\dotsm \binom{x_N}{n_N} \mu(\Delta^{\nn}f(\underline{0}))\quad\text{for all
$\underline{x} \in \Z^N$\!.}
\end{equation*}
\end{itemize}
\end{cor}
\begin{proof}
The map $\mu_*: B^{\Z^N} \!\ra (B')^{\Z^N}$ is a homomorphism of
$\Z[\Z^N]$-modules. Therefore, for all $\nn\in \N^N$\!,
\begin{equation}
\label{HOMFUNCIIEQ1}
\Delta^{\nn}\mu_* f \,=\, \Delta^{\nn}(\mu_* f) \,=\, \mu_*(\Delta^{\nn} f) \,=\, \mu_*\Delta^{\nn} f\,.
\end{equation}
From this and Theorem \ref{4.2}a) it follows that, for all $\underline{x}\in \N^N$\!,
\[ \mu_*f(\underline{x})
 \,= \sum_{\nn \in \N^N} \binom{x_1}{n_1}
 \dotsm \binom{x_N}{n_N} \Delta^{\nn}\mu_* f(\underline{0})
 \,=  \sum_{\nn \in \N^N} \binom{x_1}{n_1}
 \dotsm \binom{x_N}{n_N} \mu(\Delta^{\nn} f(\underline{0}))\,, \]
establishing part a). Applying Theorem \ref{4.2}b) to $\mu_* f$ and using
(\ref{HOMFUNCIIEQ1}) again, we get that, for all $\underline{x}\in \Z^N$\!,
\[ \mu_*f(\underline{x})
 \,= \smash[b]{\sum_{\atop{\nn \in \N^N}{\!\!|\nn| \leq d}}\binom{x_1}{n_1}
 \dotsm \binom{x_N}{n_N} \Delta^{\nn} \mu_* f(\underline{0})
 \,= \sum_{\atop{\nn \in \N^N}{\!\!|\nn| \leq d}}\binom{x_1}{n_1}
 \dotsm \binom{x_N}{n_N}  \mu( \Delta^{\nn} f(\underline{0}))\,.}\]
\end{proof}

\subsection{Polynomial Functions and Integer-Valued Polynomials}
In this section we use Theorem \ref{4.2} to compare integer-valued
polynomials to functions of finite functional degree.  The results of this section are
\emph{not} used elsewhere in this paper.  However, integer-valued polynomials and
their reductions occur in Wilson's proof of Ax-Katz over $\F_p$ \cite[Lemma
4]{Wilson06}, and the technique of representing functions between residue rings of
$\Z$ via integer-valued polynomials also occurs in a work of Varga \cite{Varga14}
generalizing Warning's Second Theorem.  It seems potentially useful to know that these
techniques can be viewed in terms of the Aichinger-Moosbauer calculus.
\\ \\
Let $R$ be a non-trivial commutative ring, let $N \in \Z^+$\!, and consider the
\textbf{evaluation map}
\[ E: R[t_1,\dotsc,t_n] \longrightarrow R^{R^N}\!, \ f \longmapsto (x \mapsto f(x)). \]
This is an $R$-algebra homomorphism; its image is, by definition, the ring of
\textbf{polynomial functions} on $R^N$\!, which we denote by $\PP(R^N\!,R)$. The
map $E$ is never an isomorphism, though the manner of the failure depends upon
$R$.   If $R$ is finite then $R[t_1,\dotsc,t_n]$ is infinite while $R^{R^N}$ is finite, so
$E$ has an infinite kernel.  If $R$ is infinite, then $E$ is not surjective \cite[Thm.\
4.3]{Clark14}. More precisely:

\begin{prop}
\label{POLYPROP1} Let $R$ be a non-trivial commutative ring, and let $N \in \Z^+$\!.
Then the following properties are equivalent:
\begin{itemize}
\item[(i)] The evaluation map $E: R[t_1,\dotsc,t_N] \longrightarrow R^{R^N}$ is
    surjective.\smallskip
\item[(ii)] The function
\[ \delta_{0,1} \in R^{R^N}\!,\ x \longmapsto \begin{cases} 1 & \text{ if } x = 0 \\ 0 & \text{ if } x
\neq 0 \end{cases} \]
lies in the image of $E$.\smallskip
\item[(iii)] The ring $R$ is a finite field.
\end{itemize}
\end{prop}
\begin{proof}
To show that (iii) implies (i), which then entails (ii), assume that $R  = \F_q$ is a finite
field. In this case the study of $E$ was the essence of Chevalley's proof of Theorem
\ref{CWAK}a) in \cite{Chevalley35}. He showed that $E$ is surjective, with kernel
$\langle t_1^q-t_1,\dotsc,t_N^q - t_N \rangle$.  For English language proofs of
modest generalizations, see
\cite[Cor.\ 2.5 and Prop.\ 4.4]{Clark14}.\smallskip  \\
To show that $\neg$(iii) implies $\neg$(ii), which then entails $\neg$(i),
we distinguish two cases.\\
\emph{Case 1, $R$ is not a field:} In this case, there exists a proper ideal $I$ in $R$,
and then every function $F$ in the image $\PP(R^N\!,R)$ of $E$ is
\textbf{congruence-preserving} module $I$: if $x = (x_1,\dotsc,x_N)$, $y =
(y_1,\dotsc,y_N) \in R^N$ are such that $x_i \equiv y_i \pmod{I}$ for all $1 \leq i \leq
N$, then $f(x) \equiv f(y) \pmod{I}$. But, if $a \in I \!\setminus\! \{0\}$, then $a \equiv 0
\pmod{I}$ while
\[ \delta_{0,1}(0,\dotsc,0) \,=\, 1 \,\not \equiv\, 0 \,=\, \delta_{0,1}(a,\dotsc,a) \pmod{I}. \]
So $\delta_{0,1}$ is not congruence-preserving and does not lie in the image of $E$.  \\
\emph{Case 2, $R$ is not finite:}  In this case, \cite[Thm.\ 4.9a)]{Clark-Schauz21}
gives $\fdeg(\delta_{0,1}) = \infty$. So by Lemma \ref{AM12.5}, $\delta_{0,1}$ is not
a polynomial function, and does not lie in the image of $E$.
\end{proof}
\noindent If $R$ is an infinite commutative ring that is not a field, we just gave two
proofs (in Cases 1 and 2) that $\delta_{0,1} \in R^{R^N} \setminus \PP(R^N\!,R)$. The
second proof showed more: that $\delta_{0,1}$ has infinite functional degree.  In
general, for non-trivial commutative rings $R$, Lemma \ref{AM12.5} says that
\begin{equation}\label{eq.F}
\PP(R^N\!,R) \,\subseteq\, \mathcal{F}(R^N\!,R) \,\subseteq\, R^{R^N}\!.
\end{equation}
This
leads to a more interesting version of the question of when $E$ is surjective.

\begin{ques}
\label{REFINEDQUES} For which non-trivial commutative rings $R$ and numbers $N
\in \Z^+$ do we have $\PP(R^N\!,R) = \mathcal{F}(R^N\!,R)$ -- i.e., when is every
function $f \in R^{R^N}\!$ of finite functional degree a polynomial function?
\end{ques}
\noindent Here is an answer to Question \ref{REFINEDQUES} when $R$ is finite.

\begin{prop}
\label{POLYPROP2} For non-trivial finite commutative rings $R$, the following
properties are equivalent:
\begin{itemize}
\item[(i)] $\PP(R^N\!,R) = \mathcal{F}(R^N\!,R)$ for all $N \in \Z^+$\!.\smallskip
\item[(ii)] $\PP(R^N\!,R) = \mathcal{F}(R^N\!,R)$ for some $N \in \Z^+$\!.\smallskip
\item[(iii)] $R \,\cong\, \prod_{i=1}^r \F_{p_i^{\alpha_i}}$ for some
    $r,\alpha_1,\dotsc,\alpha_r\in \Z^+\!$ and prime numbers $p_1 < \dotsb < p_r$.
\end{itemize}
\end{prop}
\begin{proof}
If $R$ is a finite commutative ring of order $p_1^{\alpha_1} \dotsm p_r^{\alpha_r}$
(for primes $p_1 < \dotsb < p_r$), then we have a unique internal direct product
decomposition $R = \prod_{i=1}^r \mathfrak{r}_i$\,, with $\mathfrak{r}_i$ a ring of
order $p_i^{\alpha_i}$ \cite[Thm.\ 8.37]{Clark-CA} -- the \textbf{$p_i$-primary
component} of $R$.  We have a natural ring isomorphism
\[ R[t_1,\dotsc,t_n] \,=\, \prod_{i=1}^r \mathfrak{r}_i[t_1,\dotsc,t_n] \]
and also, by \cite[Thm.\ 9.4]{Aichinger-Moosbauer21} or \cite[Thm.\
3.13]{Clark-Schauz21}, a natural decomposition
\[ \mathcal{F}(R^N\!,R) \,=\, \prod_{i=1}^r \mathcal{F}(\mathfrak{r}_i^N\!,\mathfrak{r}_i). \]
Using these decompositions we get that
\[ \PP(R^N\!,R) = \mathcal{F}(R^N\!,R)\ \iff \ \forall 1 \leq i \leq N, \ \PP(\rr_i^N,\rr_i)
 = \mathcal{F}(\rr_i^N,\rr_i). \]
So we reduce to the case in which $R$ has prime power order and, by \cite[Thm.\
9.1]{Aichinger-Moosbauer21},
\[ R^{R^N} =\, \mathcal{F}(R^N\!,R). \]
Hence, our problem reduces to the previous problem of when the evaluation map is
surjective.  By Proposition \ref{POLYPROP1}, this holds if and only if $R$ is a finite
field.  So, independent of $N \in \Z^+$: $\PP(R^N\!,R) = \mathcal{F}(R^N\!,R)$ if and
only if, for all $1\leq i\leq r$, $\mathfrak{r}_i$ is a finite field $\F_{p_i^{\alpha_i}}$, i.e.\
$R \cong \prod_{i=1}^r \F_{p_i^{\alpha_i}}$.
\end{proof}

\noindent When $R$ is infinite we do not know a complete answer to Question
\ref{REFINEDQUES}, but we will exhibit some positive and negative results.

\begin{lemma} \mbox{}
\label{4.2lem}
Let $h \in \mathcal{F}(\Q^N\!,\Q)$. If $h|_{\Z^N} = 0$ then $h = 0$.
\end{lemma}

\begin{proof}
Let $D \in \Z^+$\!, and define $h_D \in \Q^{\Z^N}$ by
\[ h_D(\underline{x}) \,\coloneqq\, h\left(\frac{x_1}{D},\dotsc,\frac{x_N}{D}\right). \]
The function $h_D$ is obtained by precomposing $h$ with a group endomorphism of
$(\Q^N\!,+)$, so $h_D \in \mathcal{F}(\Q^N\!,\Q)$ by \cite[Thm.\
4.3]{Aichinger-Moosbauer21}. Hence, by Corollary \ref{4.2cor}, there exists a
polynomial $\hat h_D(\underline{t})\in\Q[t_1,\dotsc,t_N]$ with $\hat
h_D(\underline{x})=h_D(\underline{x})$ for all $\underline{x} =
(x_1,\dotsc,x_N)\in\Z^N$\!.  Applying Lemma \ref{CATS} to $\hat{h}_D$ with $X =
(D\Z)^N$ gives $\hat{h}_D = 0$.  Thus for all $D \in \Z^+$ we have $h|_{(D^{-1}\Z)^N}
= 0$, so $h = 0$.
\end{proof}

\begin{prop}
\label{POLYPROP3} For all $N \in \Z^+$\!, we have $\PP(\Q^N\!,\Q) =
\mathcal{F}(\Q^N\!,\Q)$.
\end{prop}
\begin{proof}
That $\PP(\Q^N\!,\Q)\subseteq\mathcal{F}(\Q^N\!,\Q)$ is clear. To show that
$\mathcal{F}(\Q^N\!,\Q)\subseteq\PP(\Q^N\!,\Q)$ let $g \in \mathcal{F}(\Q^N\!,\Q)$,
say with $\fdeg(g)\leq d\in\N$. By Corollary \ref{4.2cor}, there exists a polynomial
$\hat g \in \Q[t_1,\dotsc,t_N]$ with $\deg(\hat
g)\leq d$ such that $\hat g(\underline{x})=g(\underline{x})
$, for all $\underline{x} = (x_1,\dotsc,x_N)\in\Z^N$\!. So $h:=E(\hat g)-g$ is zero on
$\Z^N$\! and $\fdeg(h)\leq d$ by \cite[Lemma 3.2]{Aichinger-Moosbauer21}. By
Lemma \ref{4.2lem} it follows that $h=0$, which implies $E(\hat g) = g$, i.e.\
$g\in\PP(\Q^N\!,\Q)$.
\end{proof}

\noindent From now until the end of \S 2.3 we will assume that $R$ is a commutative
domain of characteristic $0$, with fraction field $K$. In this case the evaluation map
$E: R[t_1,\dotsc,t_N] \longrightarrow R^{R^N}$ is injective \cite[Prop. 4.5]{Clark14}
and thus induces an isomorphism $R[t_1,\dotsc,t_N] \stackrel{\sim}{\longrightarrow}
\PP(R^N\!,R)$.  It is a result of Aichinger-Moosbauer \cite[Lemma
10.4]{Aichinger-Moosbauer21} that for all $f \in K[t_1,\dotsc,t_n]$ we have
$\fdeg(E(f)) = \deg(f)$.  We will show that the same conclusion holds over $R$ and, in
fact, a little more. Namely, we consider the subring of \textbf{integer-valued
polynomials}
\[ \Int(R^N\!,R) \,\coloneqq\, \bigl\{f \in K[t_1,\dotsc,t_N]
\mid E(f)(R^N)\subseteq R \bigr\} \,\subseteq\, K[t_1,\dotsc,t_N]. \]

\begin{prop}
\label{POLYPROP4} Let $R$ be a commutative domain. If $f \in \Int(R^N\!,R)$ and
$E_R(f):=\bigl(x \mapsto f(x)\bigr)\in R^{R^N}$\!\!, then
\begin{equation}
\label{eq.E_R}
\fdeg(E_R(f)) \,\leq\, \deg(f),
\end{equation}
with equality if $R$ has characteristic $0$.
\end{prop}

\begin{proof}
Let $f \in \Int(R^N\!,R)$. By Lemma \ref{HOMFUNCI}, domain restriction and
codomain restriction does not increase the functional degree, so Lemma \ref{AM12.5}
yields
\[ \fdeg(E_R(f)) \,\leq\, \fdeg(E(f)) \,\leq\, \deg(f). \]
Now, assume that $R$ has characteristic $0$ and that $d:=\deg(f)\geq0$. To
complete the proof, it suffices to show that $\fdeg(E_R(f))\geq d$. As $\deg(f)=d$, a
monomial $t_1^{n_1}t_2^{n_2}\dotsm t_N^{n_N}$ with $n_1+n_2+\dotsb+n_N=d$
occurs in the standard expansion of $f$. If the operators $\Delta_i$ are applied to
polynomials in the same way as they are applied to functions, then
$\deg(\Delta_1^{n_1}\dotsm \Delta_N^{n_N}f)=\deg(f)-(n_1+\dotsb+n_N)=0$,
because each application of a $\Delta_i$ reduces the degree by exactly $1$, as the
quotient field of $R$ has characteristic $0$. This shows that the function
$\Delta_1^{n_1}\dotsm \Delta_N^{n_N}E_R(f)=E_R(\Delta_1^{n_1}\dotsm
\Delta_N^{n_N}f)$ is constant but not zero, so that $\fdeg(E_R(f))\geq d$, indeed.
\end{proof}

\noindent If $R$ is a commutative domain of characteristic $p > 0$, strict inequality
can occur in \eqref{eq.E_R}. To get an equality one needs to use the $p$-weight
degree and, when $R$ is finite, reduced polynomials: see \S 4.3.
\\

\noindent For a commutative domains $R$, Proposition \ref{POLYPROP4} gives a
refinement of \eqref{eq.F}:
\begin{equation}
\PP(R^N\!,R) \,\subseteq\, \Int(R^N\!,R) \,\subseteq\, \mathcal{F}(R^N\!,R)
 \,\subseteq\, R^{R^N}\!.
\end{equation}
\noindent
This yields a negative answer to Question \ref{REFINEDQUES} whenever
$\Int(R^N\!,R) \supsetneq \PP(R^N\!,R)$, which certainly holds for $R = \Z$ as e.g.\
$t(t-1)/2$ is an integer-valued polynomial that does not lie in $\Z[t]$. This leads us to
the following result:

\begin{thm}
\label{4.3} \leavevmode
\begin{itemize}
\item[a)] $\,\mathcal{B} \coloneqq \bigl\{ \binom{x_1}{n_1} \dotsm \binom{x_n }{n_N }
    \mid \nn \in \N^N \bigr\}$ is a basis of the $\Z$-module
    $\mathcal{F}(\Z^N\!,\Z)$.\smallskip
\item[b)] $\,\mathcal{F}(\Z^N\!,\Z) = \Int(\Z^N\!,\Z)$.
\end{itemize}
\end{thm}
\begin{proof}
Part b) of Theorem \ref{4.2} implies that $\mathcal{B}$ spans $\mathcal{F}(\Z^N\!,\Z)$
as a $\Z$-module, and part a) of that theorem states the uniqueness property that
characterizes a basis.
\smallskip \\
By Proposition \ref{POLYPROP4}, we also have $\Int(\Z^N\!,\Z) \subseteq
\mathcal{F}(\Z^N\!,\Z)$, and it remains to show that $\mathcal{F}(\Z^N\!,\Z) \subseteq
\Int(\Z^N\!,\Z)$.  The well-known fact that for all $n \in \N$ we have $\binom{x}{n} \in
\Int(\Z,\Z)$ follows from Lemma \ref{4.0} and induction. Since $\Int(\Z^N\!,\Z)$ is a
ring, we have $b \in \Int(\Z^N\!,\Z)$ for all $b \in \mathcal{B}$. So,
\[ \mathcal{F}(\Z^N\!,\Z) \,=\, \langle \mathcal{B} \rangle_{\Z} \,\subseteq\, \Int(\Z^N\!,\Z).
\qedhere \]
\end{proof}

\noindent Theorem \ref{4.3} implies that $\mathcal{B}$ is a $\Z$-basis of the ring
$\Int(\Z^N\!,\Z)$ of integer-valued polynomials, a result of Ostrowski
\cite{Ostrowski19}. See \cite[Ch.\ 11]{Cahen-Chabert} for a general treatment of
$\Int(R^N\!,R)$ for commutative domains $R$.  Cahen-Chabert also address when
$\Int(R^N\!,R) = \PP(R^N\!,R)$ in \cite[\S I.3]{Cahen-Chabert}, showing in particular
that equality holds when every residue field of $R$ is infinite \cite[Cor.\
I.3.7]{Cahen-Chabert}, so e.g.\ when $R$ is a $\Q$-algebra.  Our next result implies
that, for each $N \in \Z^+$\!, $\Int(R^N\!,R) \subsetneq \mathcal{F}(R^N\!,R)$
whenever $R \supsetneq \Q$ is a $\Q$-algebra.
\\ \\
Let us say that a ring $R$ is a \textbf{Cayley ring} if the Cayley
homomorphism
\[ \mathfrak{C}: R \longrightarrow \End (R,+), \ r \longmapsto r\bullet: x \mapsto rx \]
is an isomorphism (equivalently, is surjective).

\begin{example}
\label{POLYPROP5}\leavevmode
\begin{itemize}
\item[a)] The following rings are Cayley ring:
\begin{itemize}
\item[1.] prime fields, i.e.\ $\Q$ and the finite fields $\F_p$ with $p\in\mathcal{P}$;
\item[2.] subrings of $\Q$, i.e.\ localizations of $\Z$, including $\Z$ and $\Q$.
\smallskip
\end{itemize}
\item[b)]  A commutative ring is \emph{not} Cayley if it is free of rank greater than
    $1$ as a module over some proper subring.  Thus, none of the following rings are
    Cayley rings:
\begin{itemize}
\item[1.] non-prime fields, i.e.\ fields other than $\Q$ and $\F_p$, for all
    $p\in\mathcal{P}$;
\item[2.] algebras $R$ over any field $F$ such that $F \subsetneq R$;
\item[3.] rings of integers $\Z_K$ of number fields $K \supsetneq \Q$;
\item[4.] valuation rings of $p$-adic fields $K \supsetneq \Q_p$, for any
    $p\in\mathcal{P}$.
\end{itemize}
\end{itemize}
\end{example}

\begin{prop}
\label{POLYPROP6} Let $R$ be a commutative domain of characteristic $0$.  If for
some $N \in \Z^+$ we have $\Int(R^N\!,R) = \mathcal{F}(R^N\!,R)$, then $R$ is a
Cayley ring.
\end{prop}
\begin{proof}
Proceeding by contraposition, suppose that $R$ is not a Cayley ring: this means
precisely that there is a $\Z$-linear map $L: (R,+) \ra (R,+)$ that is not of the form
$E(f)$ for a linear polynomial $f \in R[t]$.  If $K$ is the fraction field of $R$, then
moreover $L$ is not of the form $E(f)$ for a linear polynomial $f \in K[t]$: if $f = ax+b$
with $a,b \in K$, then evaluating at $0$ gives $b = 0$ and evaluating at $1$ gives $a
= L(1) \in R$.  Since $\fdeg(L) = 1$, by Proposition \ref{POLYPROP4}. $L$ is
therefore not given by any integer-valued polynomial.  This establishes the result for
$N = 1$. For each $N \in \Z^+$\!, the function $L_N: R^N \!\ra R$ with
$L_N(x_1,\dotsc,x_N) = L(x_1)$ is again $\Z$-linear, but it is not the restriction to
$R^N$ of any $K$-linear polynomial function. So $L_N \in \mathcal{F}(R^N\!,R)
\setminus \Int(R^N\!,R)$.
\end{proof}

\noindent Proposition \ref{POLYPROP6} and Example \ref{POLYPROP5} give lots of
examples in which $\Int(R^N\!,R) \subsetneq \mathcal{F}(R^N\!,R)$: e.g.\ any field $K
\supsetneq \Q$.  On the other hand, using similar arguments to the ones we have
made, one can show that $\Int(R^N\!,R) = \mathcal{F}(R^N\!,R)$ for any subring $R$
of $\Q$.

\subsection{Lifting}\label{sec.lift}
Suppose that $\mu: B \ra B'$ is a surjective homomorphism of commutative groups
and $f \in \mathcal{F}(\Z^N\!,B')$. By Theorem \ref{4.2}b), there is a unique function
$a_{\bullet}: \N^N \ra B'$ that is nonzero in at most finitely many points, we say
\textbf{finitely nonzero}, such that
\[ f(\underline{x}) = \sum_{\nn \in \N^N} \binom{x_1}{n_1} \cdots \binom{x_N}{n_N} a_{\nn}
 \quad\text{for all $\underline{x} \in \Z^N$\!.} \]
By a \textbf{lift} of $a_{\bullet}$ to $B$ (through $\mu$) we will mean a finitely
nonzero function $\tilde a_{\bullet}: \N^N\!\ra B$ such that $\mu \circ \tilde a_{\bullet}
= a_{\bullet}$.  A \textbf{proper lift} is a lift $\tilde a_{\bullet}$ that moreover satisfies,
for all $\nn\in \N^N$\!,
\[\tilde a_{\nn} = 0 \ \iff\ a_{\nn} = 0\,. \]
Proper lifts always exist, and they are unique if and only if $\mu: B \ra B'$ is an
isomorphism or $a_{\bullet}$ is identically $0$. To a proper lift $\tilde a_{\bullet}$ we
attach the function $\tilde f \in B^{\Z^N}\!$ defined by
\[ \tilde f(\underline{x}) \,:= \sum_{\nn \in \N^N}
 \binom{x_1}{n_1} \dotsm \binom{x_N}{n_N} \tilde a_{\nn}\,. \]
We see that $\fdeg(\tilde f)=\sup\{ |\nn| \mid \tilde a_{\nn} \neq 0\}=\sup \{ |\nn| \mid
a_{\nn} \neq 0\}=\fdeg(f)$. Combined with Corollary \ref{HOMFUNCII}, we have
\[ \mu\circ \tilde{f} = f \quad\text{and}\quad \fdeg(\tilde f) \,=\, \fdeg(f)\,. \]
We also call $\tilde{f}$ a \textbf{proper lift} of $f$, and may conversely say that $f$ is
the \textbf{reduction} of $\tilde{f}$.\smallskip

Moreover, the concept of a proper lift can be strengthened if $B$ and $B'$ are direct
products (or direct sums) of families $(B_j)_{j\in J}$ and $(B_j')_{j\in J}$, respectively,
and if $\mu$ is the Cartesian product (or direct product) of homomorphisms $\mu_j:
B_j\ra B'_j$\,. In that situation, if $\tilde\pi_j:B\ra B_j$ and $\pi_j:B'\!\ra B_j'$ denote
the corresponding projections, we have $\pi_j\circ\mu=\mu_j\circ\tilde\pi_j$ for each
$j\in J$. We then say $\tilde a_{\bullet}: \N^N\!\ra B$ is \textbf{coordinate-wise proper}
if each $\tilde \pi_j\circ \tilde a_{\bullet}: \N^N\!\ra B_j$ is a proper lift of $\pi_j\circ
a_{\bullet}: \N^N\!\ra B_j'$\,. Accordingly, we speak of a coordinate-wise proper lift
$\tilde f \in B^{\Z^N}\!$ of $f\in \mathcal{F}(\Z^N\!,B')$ if each $\tilde\pi_j\circ\tilde f$ is
a proper lift of $\pi_j\circ f$. This kind of proper lifts always exist, as well, as one can
construct the coordinate-wise lifts $\tilde \pi_j\circ \tilde a_{\bullet}$ first, and then
combine them into one lift $\tilde a_{\bullet}$.
\\ \\
Combining this discussion with Theorems \ref{4.2} and \ref{4.3}, we find that for each
$N,m \in \Z^+$\!, every $f \in \mathcal{F}(\Z^N\!,\,\Z/m\Z)$ is the reduction of an
integer-valued polynomial of degree $\fdeg(f)$.  In particular, this applies when for
some $p \in \mathcal{P}$ we have $m = p^\beta$ and $f$ is
$(p^{\alpha_1}\!,\dotsc,p^{\alpha_N})$-periodic for some $\alpha_1,\dotsc,\alpha_N
\in \Z^+$\!, i.e., if $f$ lies in the image of the natural map
$(\Z/p^\beta\Z)^{\bigoplus_{i=1}^N \Z/p^{\alpha_i} \Z} \ra (\Z/p^\beta\Z)^{\Z^N}$\!\!, a
situation that we are about to examine in more detail.

\begin{remark}
The fact that functions $\Z/p^\alpha \Z \ra \Z/p^\beta \Z$ can (after pullback via $\varepsilon: \Z \ra \Z/p^{\alpha} \Z$) be represented by
reductions of integer-valued polynomials is applied in work of Varga \cite{Varga14}.
In \cite{Clark-Watson18} this work was generalized to maps of the form
$\Z_K/\pp^\alpha \ra \Z_K/\pp^\beta$ where $K$ is a number field, $\Z_K$ is its ring
of integers, and $\pp$ is a nonzero prime ideal of $\Z_K$ (so that $\Z_K/\pp^\alpha$
and $\Z_K/\pp^\beta$ are finite rings of $p$-power order for some $p \in
\mathcal{P}$). Perhaps these works could be refined using considerations from the
present paper and \cite{Clark-Schauz21}.
\end{remark}

\subsection{Representation of Functions Between Finite Commutative $p$-Groups}
If $A$ is a finitely generated commutative group, then for some $N \in \Z^+$ we have
a surjective group homomorphism $\varepsilon: \Z^N \!\ra A$.  Indeed, up to a
harmless isomorphism, we may write $A$ as $\bigoplus_{i=1}^N \Z/a_i \Z$ with
parameters $1\neq a_i \in \N$ and then take
\[ \varepsilon:\, \Z^N\! \ra \bigoplus_{i=1}^N \Z/a_i \Z\,,\quad(x_1,\dotsc,x_N)
 \mapsto (x_1+a_1\Z,\dotsc,x_N+a_N\Z)\,. \]
As recalled in Lemma \ref{HOMFUNCI}, the pullback map $\varepsilon^*$ restricts to
an injective group homomorphism
\[ \varepsilon^*:\, \mathcal{F}(A,B) \hookrightarrow \mathcal{F}(\Z^N\!,B)\,.\]
and thus every $f \in \mathcal{F}(A,B)$ has the same functional degree as its
pullback to $\Z^N$\!, which by Theorem \ref{4.2}b) has a canonical series
representation.
\\ \\
For commutative groups $A$ and $B$, we recall the quantity
\[ \delta(A,B) \,\coloneqq\, \sup \{ \fdeg(f) \mid f \in B^A \}, \]
introduced in \cite{Aichinger-Moosbauer21} and further studied in
\cite{Clark-Schauz21}. It depends only on the isomorphism type of $A$ and $B$, and
moreover, as shown in \cite[Cor.\ 4.3]{Clark-Schauz21},
\[ \delta(A,B) \,=\, \delta(A,\Z/e(B) \Z). \]
When both $A$ and $B$ are nontrivial and finite, \cite[Thm.\ 4.9]{Clark-Schauz21}
says that $\delta(A,B) < \infty$ if and only if $A$ and $B$ are $p$-groups for the
same $p \in \mathcal{P}$.  Moreover, by \cite[Thm.\ 4.9c)]{Clark-Schauz21}, if $p \in
\mathcal{P}$ and $N,\beta,\alpha_1,\dotsc,\alpha_N \in \Z^+$\!, then
\begin{equation*}
 \delta\bigl(\bigoplus_{i=1}^N \Z/p^{\alpha_i} \Z, \Z/p^\beta\Z\bigr)\,=\,
 \delta_p\bigl(\underline{\alpha},\beta)\,,
\end{equation*}
where
\begin{equation*}
\label{CS4.9CEQ}
\delta_p\bigl(\underline{\alpha},\beta) \,:=\,
 \sum_{i=1}^N(p^{\alpha_i}\!-1) + (\beta-1)(p-1)p^{\max\{\alpha_1,\dotsc,\alpha_N\}-1}.
\end{equation*}

\begin{thm}
\label{THM4.57} Let $p \in \mathcal{P}$, let $N,\alpha_1,\dotsc,\alpha_N\in \Z^+$\!,
and put $A \coloneqq \bigoplus_{i=1}^N \Z/p^{\alpha_i} \Z$.  Let $B$ be a
commutative group, and let $F: \Z^N\! \ra B$ be the pullback of a function $f: A \ra
B$.
\begin{itemize}
\item[a)] If $\beta \in \Z^+$\! is such that $p^\beta f(a) = 0$ for all $a \in A$, then
\[F(\underline{x})\ =\!\!\! \sum_{\atop{\nn \in \N^N}{\ |\nn|
 \leq \delta_p(\underline{\alpha},\beta)\!\!\!}}
 \binom{x_1}{n_1} \dotsm \binom{x_N}{n_N} \Delta^{\nn} F(\underline{0})
 \quad\text{for all $\underline{x}\in \Z^N$\!.} \]
\item[b)] For all $h \in \Z^+$\! and all $\underline{n} \in \N^n$ with $|\underline{n}| >
    \delta_p(\underline{\alpha},h)$,
\[ \Delta^{\nn}F(\underline{x}) \,\in\, p^{h} B
 \quad\text{for all $\underline{x}\in \Z^N$\!.} \]
\item[c)] Let $\mu_h: B \ra B/p^h B$ be the quotient map.  The conclusion of part b)
    continues to hold for every function $F: \Z^N\!\ra B$ such that $\mu_h\circ F:
    \Z^N\!\ra B/p^h B$ is the pullback of a function $g: A \ra B/p^h B$.
\end{itemize}
\end{thm}
\begin{proof} a)
Let $\underline{B} \coloneqq \langle f(A) \rangle$ be the subgroup generated by the
image of $f$.  Because $f(a) \in B[p^{\beta}]$ for all $a \in A$, we have that
$\underline{B} = \underline{B}[p^{\beta}]$.  We may view $f$ as a function with
codomain $\underline{B}$, which by \cite[Cor.\ 3.10b)]{Clark-Schauz21} does not
change its functional degree, so we may assume that $B=\underline{B}$.  So by
\cite[Thm.\ 4.9c)]{Clark-Schauz21}, we get
\[\fdeg(F) \,=\, \fdeg(f) \,\leq\,\delta_p(\underline{\alpha},\beta), \]
and the result follows from Theorem \ref{4.2}.\smallskip  \\
b) Let $\mu_h: B \ra B/p^h B$ be the quotient map.  The map $\mu_h\circ f: A \ra
B/p^h B$ has functional degree at most $\delta_p(\underline{\alpha},h)$, hence so
does its pullback to $\Z^N$\!, which is $\mu_h\circ F$.
For all $\nn \in \N^N$ with $|\nn| > \delta_p(\underline{\alpha},h)$ this means
$\mu_h\circ\Delta^{\nn}F = \Delta^{\nn}(\mu_h\circ F) = 0$, which implies that, for all
$\underline{x}\in \Z^N$\!, $\mu_h(\Delta^{\nn}F(\underline{x})) = \mu_h\circ
\Delta^{\nn}F(\underline{x}) = 0$,
i.e.\ $\Delta^{\nn} F(\underline{x}) \in p^h B$.\smallskip \\
c) The proof of the previous part used only that $\mu_h\circ F$ is pulled back from
$A$.
\end{proof}
\noindent
We deduce following results, the latter being a vector-valued analogue of the former.

\begin{cor}
\label{COR4.57}  Let $p \in \mathcal{P}$, and let $N,\beta,\alpha_1,\dotsc,\alpha_N\in
\Z^+$\!.
Let  $f: \bigoplus_{i=1}^N \Z/p^{\alpha_i} \Z \ra \Z/p^{\beta}\Z$ be any function, let $F:
\Z^N\! \ra \Z/p^{\beta} \Z$ be the pullback of $f$, and let $\tilde{F}: \Z^N\! \ra \Z$ be a
proper lift of $F$. Then:
\begin{itemize}
\item[a)]\ \vspace{-.5em}
\begin{equation*}
\tilde F(\underline{x})\
=\!\!\! \sum_{\atop{\nn \in \N^N}{\ \ |\nn| \leq \delta_p(\underline{\alpha},\beta)}}\!\!\!\!
\binom{x_1}{n_1} \dotsm \binom{x_N}{n_N} \Delta^{\nn} \tilde F(\underline{0})
\quad\text{for all $\underline{x} \in \Z^N$\!.}
\end{equation*}
\item[b)] For all $h \in \Z^+\!$ and all $\underline{n} \in \N^n$ with $|\underline{n}| >
    \delta_p(\underline{\alpha},h)$,
\[ p^h \bigm{|} \Delta^{\nn}\tilde F(\underline{0})\,. \]
\end{itemize}
\end{cor}
\begin{proof}
a) Since $\fdeg(\tilde{F}) = \fdeg(F) = \fdeg(f) \leq \delta_p(\underline{\alpha},\beta)$,
this follows from Theorem \ref{4.2}b).\smallskip \\
b) Assuming $\nn>\delta_p(\underline{\alpha},h)$, we prove that $p^h$ divides
$\Delta^{\nn}\tilde F(\underline{0})$. If $h \geq \beta$ then
$$\fdeg(\tilde F) \,=\, \fdeg(F) \,=\, \fdeg(f) \,\leq\, \delta_p(\underline{\alpha},\beta)
 \,\leq\, \delta_p(\underline{\alpha},h) \,<\, |\nn|\,,$$
so that $\Delta^{\nn} \tilde F(\underline{0}) = 0$, which is divisible by $p^h$\!. Hence,
we may assume $1 \leq h < \beta$. We show that, in this case, Theorem
\ref{THM4.57}c) applies to $\tilde F$ and $\Z$ in the place of $F$ and $B$, which
then yields $\Delta^{\nn}\tilde F(\underline{0}) \,\in\, p^{h}\Z$\,, i.e.\ $p^h \mid
\Delta^{\nn}\tilde F(\underline{0})$, as desired. With the canonic surjections
$\mu_\beta: \Z \ra \Z/p^\beta \Z$ and $\mu^\beta_h: \Z/p^\beta\Z \ra \Z/p^h \Z$, it
suffices to recognize $\mu_h\circ \tilde F$ as the pullback of $g:=\mu^\beta_h\circ f$.
Since $F$ is the pullback of $f$, however, we obtain $\mu^\beta_h\circ F$ as the
pullback of $\mu^\beta_h\circ f$. But, $\mu^\beta_h\circ F=\mu_h\circ \tilde F$,
because $\mu_h = \mu^\beta_h\circ\mu_\beta$ and $\mu_\beta\circ\tilde F = F$, as
$\tilde F$ is a lift of $F$. So, indeed, $\mu_h\circ \tilde F$ is the pullback of the
function $g$.
\end{proof}

\begin{cor}
\label{COR4.57b}  Let $p \in \mathcal{P}$, let $N, \beta, \alpha_1,\dotsc,\alpha_N\in
\Z^+$\!. For each $j$ in a nonempty index set $J$, let $\beta_j\in\{1,\dotsc,\beta\}$.
Let  $f: \bigoplus_{i=1}^N \Z/p^{\alpha_i} \Z \ra \prod_{j\in J} \Z/p^{\beta_j}\Z$ be any
function, let $F: \Z^N\! \ra \prod_{j\in J} \Z/p^{\beta_j}\Z$ be the pullback of $f$, and
let $\tilde F: \Z^N\! \ra \prod_{j\in J} \Z=\Z^J$ be a proper lift of $F$. Then:
\begin{itemize}
\item[a)]\ \vspace{-.5em}
\begin{equation*}
\tilde F(\underline{x})\
=\!\!\! \sum_{\atop{\nn \in \N^N}{\ \ |\nn| \leq \delta_p(\underline{\alpha},\beta)}}\!\!\!\!
\binom{x_1}{n_1} \dotsm \binom{x_N}{n_N} \Delta^{\nn} \tilde F(\underline{0})
\quad\text{for all $\underline{x} \in \Z^N$\!.}
\end{equation*}
\item[b)] If $\tilde F$ is coordinate-wise proper
    then, for all $h\in\Z^+\!$ 
    and all $\underline{n} 
    $ with $|\underline{n}| > \delta_p(\underline{\alpha},h)$, $p^h$ divides each
    coordinate of $\Delta^{\nn}\tilde F(\underline{0})\in\Z^J$\!.
\end{itemize}
\end{cor}
\begin{proof}
a) The exponent of $\prod_{j \in J} \Z/p^{\beta_j} \Z$ divides $p^{\beta}$, so that
\[ \fdeg(\tilde{F}) \,=\, \fdeg(F) \,=\, \fdeg(f) \,\leq\, \delta_p(\underline{\alpha},\beta)\,, \]
and Theorem \ref{4.2}b) applies to give the result.\smallskip \\
b) 
Fix $j\in J$ and assume $\tilde F$ as in the hypothesis. Then $\tilde\pi_j\circ\tilde F$
is a proper lift of $\pi_j\circ F$, where $\pi_j: \prod_{i\in J}
\Z/p^{\beta_i}\Z\ra\Z/p^{\beta_j}\Z$ and $\tilde\pi_j: \Z^J\ra\Z$ are the coordinate
projections. So, the functions $\tilde F_j := \tilde\pi_j\circ\tilde F$, $F_j:= \pi_j\circ F$,
and $f_j := \pi_j \circ f$ meet the requirements of Corollary \ref{COR4.57}b), which
yields
\begin{equation*}
p^h \,\mid\, \Delta^{\nn}\tilde F_j(\underline{0})
 \,=\, \Delta^{\nn}(\tilde\pi_j\circ\tilde F)(\underline{0})
 \,=\, (\tilde\pi_j\circ\Delta^{\nn}\tilde F)(\underline{0})
 \,=\, \tilde\pi_j(\Delta^{\nn}\tilde F(\underline{0}))\,.  \qedhere
\end{equation*}
\end{proof}

\section{The Group-Theoretic Ax-Katz Theorem}

\subsection{Wilson's Lemma} Let $N \in \Z^+$\!.  For $s, t_1,\dotsc,t_N \in \N$, we put
\[ [s) \,\coloneqq\, \{0,1,\dotsc,s-1\}\quad\ \text{and}
\quad\ [s^{\underline{t}}) \,\coloneqq\, \prod_{i=1}^N [s^{t_i}). \]
\\
With $\Z_{(p)}$ we denote the set of rational numbers of non-negative $p$-adic
valuation. For each $x \in \R$, we set
\[ \overline{x} \,\coloneqq\, \max(x,0). \]
Now, let $A$ and $B$ be commutative groups, and let $S \subseteq A$ be a finite
subset. Following \cite{Karasev-Petrov12}, for each $f \in B^A$\!, we define
\[ \int_S f \,\coloneqq\, \sum_{x \in S} f(x) \in B. \]
The following result is an equivalent (but simpler) reformulation of \cite[Lemma
4]{Wilson06}.

\begin{lemma}
\label{Lem5.1} \label{5.1} Let $p \in \mathcal{P}$ and let $N,\beta \in \Z^+$\!.  If $f \in
\Z^{\Z^N}$\! is such that
\[ \fdeg(f) \,<\, (p-1)(N-\beta+1)\,, \]
then
\[ \int_{[p)^N} f 
\,\equiv\, 0 \pmod{p^\beta}. \] 
\end{lemma}
\begin{proof}
\emph{Step 1:} If $0 \leq i \leq p-2$ then $\sum_{x \in \Z/p\Z} x^i = 0$: indeed, upon
choosing a generator $\zeta$ of the cyclic group $(\Z/p\Z)^{\times}$\!, we get
\[\sum_{x \in \Z/p\Z} x^i \,=\, \sum_{j=0}^{p-2} (\zeta^i)^j \,=\,
 \frac{(\zeta^i)^{p-1}-1}{\zeta^i-1} \,=\, 0\,. \]
It follows that if $i_1,\dotsc,i_\beta\in [p-1)$ then
\[ \sum_{(x_1,\dotsc,x_\beta) \in [p)^\beta}\!\! x_1^{i_1} \dotsm x_\beta^{i_\beta}
 \,=\, \prod_{j=1}^\beta \sum_{x_j \in [p)} x_j^{i_j} \,\equiv\, 0 \pmod{p^\beta}. \]
We deduce that if $g \in \Z_{(p)}[x_1,\dotsc,x_\beta]$
has $\deg_j(g) \leq p-2$ for all $1 \leq j \leq \beta$, then
\[\int_{[p)^\beta} g \,=  \!\sum_{(x_1,\dotsc,x_\beta) \in [p)^\beta}\!\! g(x_1,\dotsc,x_\beta)
\,\equiv\, 0 \pmod{p^\beta}. \] \emph{Step 2:} If the result holds for a set of functions
$f_1,\dotsc,f_m$ then it holds for the $\Z$-submodule of $\Z^{\Z^N}\!$ that they
generate.  Because of this and Theorem \ref{4.2}, it suffices to show that the result
holds for the polynomial
\[ f \,:=\,\binom{x_1}{n_1} \dotsm \binom{x_N}{n_N}\,\in\,\Q[x_1,\dotsc,x_N]\,, \]
where $(n_1,\dotsc,n_N) \in \N^N$ is arbitrary with $|\nn| < (p-1)(N-\beta+1)$.
Under that requirement we have
\[ \# \bigl\{1 \leq j \leq N \mid n_j < p-1 \bigr\} \,\geq\, \beta\,, \]
for if not, we would get $|\nn|=\sum_{j=1}^Nn_j\geq (p-1)(N-\beta+1)$.  So, we may
assume, without loss of generality, that $n_j < p-1$ for all $1 \leq j \leq \beta$. Then,
for every fixed $\underline{y}=(y_{\beta+1},\dotsc,y_N )\in \Z^{N-\beta}$,
\[ g_{\underline{y}}(x_1,\dotsc,x_\beta)
 \,:=\,f(x_1,\dotsc,x_\beta,y_{\beta+1},\dotsc,y_N) \,\in\, \Z_{(p)}[x_1,\dotsc,x_\beta], \]
and $\deg_j(g_{\underline{y}})=n_j \leq p-2$ for all $1 \leq j \leq \beta$. So, using Step
1, we get
\[ \int_{[p)^N} f \,\,=\!\! 
\sum_{\underline{y} \in [p)^{N-\beta}}\,
\sum_{(x_1,\dotsc,x_\beta) \in [p)^\beta}\!\! g_{\underline{y}}(x_1,\dotsc,x_\beta)
\,\,\equiv\,\, 0 \pmod{p^\beta}. \qedhere\]
\end{proof}

\subsection{The Proof of Theorem \ref{AXKATZWILSON}}

\begin{proof}[Proof of Theorem \ref{AXKATZWILSON}]
Without lose of generality, we may assume that, for $1\leq j\leq r$,
\[ d_j \,\coloneqq\, \fdeg(f_j) \,>\,0\,. \]
We set
\[ \mathcal{M} \,\coloneqq\, \max_{1 \leq j \leq r} p^{\beta_j-1} d_j\,\in\, \Z^+ \]
and put
\[ \beta \,\coloneqq\, \biggl{\lceil}\frac{N- \sum_{j=1}^r \frac{p^{\beta_j}-1}{p-1}
d_j}{\max\limits_{1 \leq j \leq r} p^{\beta_j-1} d_j} \biggr{\rceil} \,=\,
\biggl{\lceil} \frac{N- \sum_{j=1}^r \frac{p^{\beta_j}-1}{p-1} d_j}{\mathcal{M}}
\biggr{\rceil}. \]
We have
\[ \beta \,<\, \frac{N- \sum_{j=1}^r \frac{p^{\beta_j}-1}{p-1} d_j}{\mathcal{M}} + 1 \]
and thus
\begin{equation}
\label{MAINPROOFEQ1}
 \sum_{j=1}^r \frac{p^{\beta_j}-1}{p-1} d_j\,<\,  N - \mathcal{M}(\beta-1)\,.
\end{equation}
For $1 \leq j \leq r$, we define $\chi_j: \Z \ra
\Z/p^{\beta} \Z$ by
\[ \chi_j(x) \,:=\, \begin{cases}
1 & \text{if }x \equiv 0 \pmod{p^{\beta_j}}, \\ 0 & \text{otherwise.} \end{cases}\] Since
$\chi_j$ is pulled back from $\Z/p^{\beta_j} \Z$, it has finite functional degree; let
$\widetilde{\chi}_j$ be a proper lift of $\chi_j$ from $\Z/p^{\beta} \Z$ to $\Z$\,. Let
$\chi: \Z^r \!\ra \Z/p^{\beta} \Z$ be the tensor product $\bigotimes_{j=1}^r \chi_j$ of the
$\chi_j$, and let $\widetilde{\chi}: \Z^r \!\ra \Z$ be the tensor product
$\bigotimes_{j=1}^r \widetilde{\chi}_j$ of the $\widetilde{\chi}_j$: for all
$(x_1,\dotsc,x_r)\in\Z^r$\!,
\[ \chi(x_1,\dotsc,x_r) \,:=\, \prod_{j=1}^r \chi_j(x_j)
\,=\, \prod_{j=1}^r \widetilde{\chi}_j(x_j) + p^{\beta} \Z \, =\,
 \widetilde{\chi}(x_1,\dotsc,x_r) + p^{\beta} \Z\,. \]
For $1 \leq j \leq r$, let $F_j: \Z^N \!\ra \Z/p^{\beta_j} \Z$ be the pullback of $f_j$ and
let $\tilde{F_j}$ be a proper lift of $F_j$ from $\Z/p^{\beta_j} \Z$ to $\Z$.  Then, with
the bijection $[p)^N\!\to (\Z/p\Z)^N$ given by
$\underline{x}\mapsto[\underline{x}]:=(x_j+p\Z)_{j=1}^N$, we have for all
$\underline{x}\in [p)^N$\!,
\[ \chi(\tilde{F}_1(\underline{x}),\dotsc,\tilde{F}_r(\underline{x}))
\,=\, \begin{cases} 1 & \text{if}\ [\underline{x}] \in Z(f_1,\dotsc,f_r), \\ 0 &
\text{otherwise.} \end{cases} \] Thus, the desired conclusion that $p^\beta$
divides $\# Z(f_1,\dotsc,f_r)$ is equivalent to
\[ \int_{[p)^N} \chi(\tilde{F}_1,\dotsc,\tilde{F}_r) \,=\, 0 \,\in\, \Z/p^\beta \Z\,, \]
and thus also to
\begin{equation}
\label{MAINPROOFEQ2}
 \ord_p \biggl(\int_{[p)^N} \widetilde{\chi}(\tilde{F}_1,\dotsc,\tilde{F}_r)\biggr)
 \,\geq\, \beta\,.
\end{equation}
By Corollary \ref{COR4.57}, there is a function $c_j: \N \ra \Z$ with $c_j(n) = 0$ for all
but finitely many $n\in\N$, such that, for all $x \in \Z$,
\[ \tilde{\chi}_j(x) \,=\, \sum_{n \in \N} \binom{x_j}{n} c_j(n)\,, \]
and, for each $h \in \Z^+$\!,
\begin{equation}
\label{MAINPROOFEQ3}
\ n_j > p^{\beta_j}-1 + (h-1)(p-1) p^{\beta_j-1}
\,\implies\, p^h \bigm{|} c_j(n_j)\,.
\end{equation}
Therefore, for all $\underline{x}\in\Z^N$\!,
\begin{eqnarray*}
\tilde{\chi}(\tilde{F}_1(\underline{x}),\dotsc,\tilde{F}_r(\underline{x})) &=&
\tilde{\chi}_1(\tilde{F_1}(\underline{x})) \dotsm \tilde{\chi}_r(\tilde{F_r}(\underline{x})) \\
&=&\!\sum_{\nn \in \N^r} \binom{ \tilde{F}_1(\underline{x})}{n_1} \dotsm
\binom{\tilde{F}_r(\underline{x})}{n_r}c_1(n_1) \dotsm c_r(n_r)\,.
\end{eqnarray*}
Hence,
\begin{eqnarray*}
\int_{[p)^N} \tilde{\chi}(\tilde{F}_1,\dotsc,\tilde{F}_r) &=& \int_{[p)^N}
\sum_{\nn \in \N^r} \binom{ \tilde{F}_1}{n_1} \dotsm
\binom{\tilde{F}_r}{n_r}c_1(n_1) \dotsm c_r(n_r) \\
&=&\!\sum_{\nn \in \N^r} c_1(n_1) \dotsm c_r(n_r)
\int_{[p)^N}\binom{ \tilde{F}_1}{n_1} \dotsm \binom{\tilde{F}_r}{n_r}\,.
\end{eqnarray*}
Thus to prove (\ref{MAINPROOFEQ2}), it suffices to show that, for each $\nn =
(n_1,\ldots,n_r) \in \N^r$\!,
\begin{equation}
\label{MAINPROOFEQ4}
\ord_p \bigl(
c_1(n_1)\dotsm c_r(n_r)\bigr) \,+\, \ord_p \biggl( \int_{[p)^N}
\binom{ \tilde{F}_1}{n_1} \dotsm \binom{\tilde{F}_r}{n_r}\biggr)\,\,\geq\,\, \beta.
\end{equation}
Now, for each $1 \leq j \leq r$, let $h_j$ be the unique integer with
\begin{equation}
\label{HJEQ}
 p^{\beta_j}-1 + (h_j-1)(p-1) p^{\beta_j-1} \ <\  n_j \ \leq\
p^{\beta_j}-1 + h_j(p-1) p^{\beta_j-1}.
\end{equation}
Then
(\ref{MAINPROOFEQ3}) gives 
$\ord_p(c_j(n_j)) \geq h_j$. 
So, on one hand,
\begin{equation}
\label{MAINPROOFEQ5}
 \ord_p \bigl(
c_1(n_1) \dotsm c_r(n_r)\bigr)
\,=\, \sum_{j=1}^r \ord_p (c_j(n_j)) \,\geq\, \sum_{j=1}^r h_j \,=:\, \alpha\,,
\end{equation}
and then (\ref{MAINPROOFEQ4}) certainly holds if $\alpha \geq \beta$. On the other
hand, if $\alpha < \beta$ then $\beta-1-\alpha\geq0$. Hence, as $\fdeg(\tilde{F}_j) =
\fdeg(F_j)= \fdeg(f_j) \leq d_j$, using  \cite[Thm.\,4.3 \&
Lem.\,6.1]{Aichinger-Moosbauer21}, (\ref{HJEQ}), the definition of $\mathcal{M}$,
(\ref{MAINPROOFEQ1}) and that $\mathcal{M} \geq 1$, we get


\begin{align*}
\fdeg\biggl( \binom{\tilde{F}_1}{n_1} \dotsm \binom{\tilde{F}_r}{n_r}\biggr)
&\,\leq\,\,\sum_{j=1}^r n_j d_j\\
&\,\leq\,\, (p-1) \sum_{j=1}^r \Bigl( \frac{p^{\beta_j}-1}{p-1} +
h_j p^{\beta_j-1} \Bigr) d_j\\
&\,\leq\,\, (p-1)\Bigl(\sum_{j=1}^r \frac{p^{\beta_j}-1}{p-1}d_j +
\mathcal{M}\alpha\Bigr)\\
&\,<\,\, (p-1)\bigl(N - \mathcal{M}(\beta-1-\alpha)\bigr)\\[2pt]
&\,\leq\,\, (p-1)\bigl(N - (\beta-\alpha)+1\bigr).
\end{align*}
Hence, Lemma \ref{5.1} implies
\begin{equation}
\label{MAINPROOFEQ6}
\ord_p\biggl( \int_{[p)^N} \binom{\tilde{F_1}}{n_1} \dotsm \binom{\tilde{F_r}}{n_r} \biggr)
\,\geq\,\, \beta-\alpha\,.
\end{equation}
Combining (\ref{MAINPROOFEQ5}) and (\ref{MAINPROOFEQ6}) we get
(\ref{MAINPROOFEQ4}), which completes the proof of Theorem
\ref{AXKATZWILSON}.\bigskip
\end{proof}

\section{$p$-weights}

\subsection{$p$-weight degrees}
Let $p \in \mathcal{P}$. Each $d \in \N$ can be written in the form $d =
\sum_{i=0}^N a_i p^i$ with uniquely determined coefficients $a_i \in [p)$.  Using this
base $p$ expansion, we define the $p$-\textbf{weight} of $d$ as
\[ \sigma_p(d) \,=\, \sigma_{p,\N}(d) \,\coloneqq\, \sum_{i=0}^N a_i. \]
We have $\sigma_p(d) \leq d$ with equality if and only if $d \in [p)$.  For fixed $p$
and large $d$, we have $\sigma_p(d) = O(\log d)$, so the $p$-weight of $d$ can be
much smaller than $d$ itself.
\\ \\
Let $R$ be a commutative rng.  The $p$-\textbf{weight degree} of a nonzero
monomial term $c\,t_1^{d_1} \dotsm t_n^{d_n}$ with $c \in R \setminus \{0\}$ is
defined to be
\[ \sigma_p(c\, t_1^{d_1} \dotsm t_n^{d_n}) \,\coloneqq\, \sum_{i=1}^n \sigma_p(d_i), \]
and the $p$-weight degree of a nonzero polynomial $f \in R[t_1,\dotsc,t_n]$ is the
maximum $p$-weight degree of its nonzero monomial terms.  We also set
$\sigma_p(0) \coloneqq -\infty$. A polynomial has positive degree if and only if it has
positive $p$-weight degree.\\ \\
We will also need the product $\bigotimes_{i=1}^n f_i$ of functions $f_1: A_1 \ra R$\,,
\dots, $f_n: A_n \ra R$ on commutative groups $A_1,\dotsc,A_n$, where $R$ is
again a rng, which is defined by
$$\bigotimes_{i=1}^n f_i:\ \prod_{i=1}^n A_i \ra R\,,
 \quad(x_1,\dotsc,x_n)\mapsto f_1(x_1) \dotsm f_n(x_n)\,.$$
In this setting, we have the following lemmas.

\begin{lemma}
\label{AM6.2a} For each $1\leq j\leq n$, let $a_{j,1},\dotsc,a_{j,K(j)}\in A_j\subseteq
\prod_{i=1}^n A_i$\,, and let $(a_1,\dotsc,a_K)$ be a permutation of all
$K:=K(1)+\dotsb+K(n)$ given elements $a_{j,k}$\,. Then
\[ \Delta_{a_{1}}\!\dotsm\Delta_{a_{K}}\bigl(f_1\otimes\dotsm\otimes f_n\bigr)
 \,=\, \bigl(\Delta_{a_{1,1}}\!\dotsm\Delta_{a_{1,K(1)}}f_1\bigr)\otimes\dotsm\otimes
 \bigl(\Delta_{a_{n,1}}\!\dotsm\Delta_{a_{n,K(n)}}f_n\bigr)\,.\]
\end{lemma}
\begin{proof}
If $a=(a,0)\in A_1=A_1\times\{0\}\subseteq A_1\times A_2$ and $(x_1,x_2)\in
A_1\times A_2$ then
\begin{align*}
\bigl(\Delta_{a}(f_1\otimes f_2)\bigr)(x_1,x_2)
 &\,=\,f_1(x_1+a) f_2(x_2) - f_1(x_1) f_2(x_2)\\
 &\,=\,(f_1(x_1+a) - f_1(x_1) )f_2(x_2)\\
 &\,=\,\bigl((\Delta_{a}f_1)\otimes f_2\bigr)(x_1,x_2)\,.
\end{align*}
Hence, $\Delta_{a}(f_1\otimes f_2)=(\Delta_{a}f_1)\otimes f_2$\,. More generally, if
$a\in A_j$ then
$$\Delta_{a}(f_1\otimes\dotsm\otimes f_n)
 \,=\,f_1\otimes\dotsm\otimes f_{j-1}\otimes(\Delta_{a}f_j)
 \otimes f_{j+1}\otimes\dotsm\otimes f_n\,.$$
From this, and the commutativity of the operators $\Delta_{a_{j,k}}$, the stated
equation follows.
\end{proof}

\begin{lemma}
\label{AM6.2} We have
\[ \fdeg_j\bigl(\bigotimes_{i=1}^n f_i\bigr) \,\leq\, \fdeg(f_j)\quad\text{for all $\,1\leq j\leq n\,$.}\]
In particular,
\[ \fdeg\bigl(\bigotimes_{i=1}^n f_i\bigr) \,\leq\, \sum_{i=1}^n \fdeg(f_i)\,.\]
Equality holds in both inequalities, as shown in \cite[Lemma
6.2]{Aichinger-Moosbauer21}, if $R$ is a domain and the functions $f_1$, \dots,
$f_n$ are all nonzero.
\end{lemma}
\begin{proof}
Assume $1\leq j\leq n$. Lemma \ref{AM6.2a} shows that
$\Delta_{a_{j,1}}\!\dotsm\Delta_{a_{j,K(j)}}\bigl(f_1\otimes\dotsm\otimes f_n\bigr)=0$
whenever $K(j)>\fdeg(f_j)$, because
$$\Delta_{a_{j,1}}1\!\dotsm\Delta_{a_{j,K(j)}}\bigl(f_1\otimes\dotsm\otimes f_n\bigr)
 \,=\, f_1\otimes\dotsm\otimes f_{j-1}\otimes\dotsm\otimes
 \bigl(\Delta_{a_{j,1}}\!\dotsm\Delta_{a_{j,K(j)}}f_j\bigr)\otimes
 f_{j+1}\otimes\dotsm\otimes f_n$$
and $\Delta_{a_{j,1}}\!\dotsm\Delta_{a_{j,K(j)}}f_j=0$ if $K(j)>\fdeg(f_j)$. This means
$\fdeg_j\bigl(\bigotimes_{i=1}^n f_i\bigr)\leq\fdeg(f_j)$. If we combine these
inequalities with Theorem \ref{thm.pfdeg}, we get
\[ \fdeg\bigl(\bigotimes_{i=1}^n f_i\bigr)
 \,\leq\, \sum_{i=1}^n \fdeg_i\bigl(\bigotimes_{j=1}^n f_j\bigr)
 \,\leq\, \sum_{i=1}^n \fdeg(f_i)\,.\]
It remains to show that
\[ \fdeg\bigl(\bigotimes_{i=1}^n f_i\bigr)
 \,\geq\, \sum_{i=1}^n \fdeg(f_i)\]
whenever $R$ is a domain and $K(j):=\fdeg(f_j)\geq0$, for all $1\leq j\leq n$. To
prove this, we choose for each $1\leq j\leq n$ elements $a_{j,1},\dotsc,a_{j,K(j)}\in
A_j$ such that
$$\Delta_{a_{j,1}}\!\dotsm\Delta_{a_{jK(j)}}f_j\neq0\,.$$
Then, by Lemma \ref{AM6.2a}, and because $R$ is a domain,
$$\Delta_{a_{1,1}}\!\dotsm\Delta_{a_{n,K(n)}}\bigl(f_1\otimes\dotsm\otimes f_n\bigr)
 \,=\,\bigl(\Delta_{a_{1,1}}\!\dotsm\Delta_{a_{1,K(1)}}f_1\bigr)\otimes\dotsm\otimes
 \bigl(\Delta_{a_{n,1}}\!\dotsm\Delta_{a_{n,K(n)}}f_n\bigr)\,\neq\,0\,,$$
which means that $\fdeg\bigl(\bigotimes_{i=1}^n f_i\bigr)\geq K(1)+\dotsb+
K(n)=\sum_{i=1}^n \fdeg(f_i)$, indeed.
\end{proof}

\noindent The next result is the first half of \cite[Theorem
10.3]{Aichinger-Moosbauer21} in a more general setting.

\begin{prop}
\label{PWEIGHTPROP} Let $p \in \mathcal{P}$, and let $R$ be a commutative ring of
characteristic $p$.  Let $f \in R[t_1,\dotsc,t_n]$ be a polynomial, with associated
function $E(f) \in R^{R^n}$\!\!.  Then
\begin{equation}
\label{PWEIGHTINEQ}
 \fdeg(E(f)) \,\leq\, \sigma_p(f)\,.
\end{equation}
\end{prop}

\begin{proof}
Since $\fdeg(E(f)) = -\infty$ if $E(f) = 0$, we may assume that $E(f) \neq 0$.  By
\cite[Lemma 3.2]{Aichinger-Moosbauer21} we have $\fdeg(f_1+f_2) \leq
\max\bigl(\fdeg(f_1),\fdeg(f_2)\bigr)$.  Since $\sigma_p(f)$ is the maximum of the
$p$-weight degrees of the nonzero monomial terms of $f$, we reduce to the case of
a monomial term
\[ f \,=\, c\, t_1^{d_1} \dotsm t_n^{d_n},\, \ c \in R \setminus \{0\}\,. \]
Using \cite[Lemmas 6.1]{Aichinger-Moosbauer21} and Lemma \ref{AM6.2}, we get
\[ \fdeg(c\, t_1^{d_1} \dotsm t_n^{d_n}) \,\leq\, \fdeg(c) + \sum_{i=1}^n \fdeg(E(t_i^{d_i}))
 \,=\, \sum_{i=1}^n \fdeg(E(t_i^{d_i}))\,. \]
We have reduced to the univariate monomial case and must show: for all $d \in \Z^+$
we have
\[ \fdeg(E(t^d)) \,\leq\, \sigma_p(d)\,. \]
Writing $d = \sum_{i=0}^N a_i p^i$ with $a_i \in [p)$ and using \cite[Lemma
6.1]{Aichinger-Moosbauer21}, we get
\[ \fdeg(E(t^d)) \,=\, \fdeg\biggl(\prod_{i=0}^N(E(t^{p^i}))^{a_i}\biggr)
 \,\leq\, \sum_{i=0}^N a_i \fdeg(E(t^{p_i})) \,=\, \sum_{i=0}^{N}a_i \,=\, \sigma_p(d)\,, \]
since each $E(t^{p^i})$ is a nonzero group homomorphism and thus has functional
degree $1$.
\end{proof}

\subsection{A Generalized Moreno-Moreno Theorem}
Combining Corollary \ref{GTPAKT} and Proposition \ref{PWEIGHTPROP} we get:

\begin{thm}
\label{GENMORENOMORENO} Let $R$ be a finite commutative ring of prime
characteristic $p$ and order $p^N$\!.  Let $f_1,\ldots,f_r \in R[t_1,\ldots,t_n]$ be
nonzero polynomials.   If $Z:=Z_{R^n}(f_1,\dotsc,f_r)$, then
\begin{equation*}
\ord_p(\# Z) \,\geq\, \bigg{\lceil} \frac{N\bigl(n-\sum_{j=1}^r\sigma_p(f_j)\bigr)}
{\max_{j=1}^r \sigma_p(f_j)} \bigg{\rceil}.
\end{equation*}
\end{thm}

\begin{proof}
The result trivially holds if all functions $E(f_j)$ are zero. If some but not all
functions $E(f_j)$ are zero, it is enough to prove the theorem for the set of functions
$f_j$ with $E(f_j)\neq0$, as that yields a lower bound at least as good as the stated
one. So we may assume that Corollary \ref{GTPAKT} applies. The resulting inequality
$$\ord_p(\# Z) \,\geq\, \biggl{\lceil} \frac{N\bigl(n-\sum_{j=1}^r\fdeg(E(f_j))\bigr)}
 {\max_{j=1}^r \fdeg(E(f_j))}\biggr{\rceil}$$
remains true if every functional degree $\fdeg(E(f_j))$ is replaced by an upper bound
for $\fdeg(E(f_j))$, such as the one given in Proposition \ref{PWEIGHTPROP}.
\end{proof}

\noindent
If in Theorem \ref{GENMORENOMORENO} we take $R$ to be the finite field $\F_{p^N}$, we recover the Moreno-Moreno Theorem
\cite[Thm.\ 1]{Moreno-Moreno95}.

\begin{thm}[Moreno-Moreno]
\label{MORENOMORENOTHM} Let $p\in\P$ and $q:=p^N$. Let $f_1,\ldots,f_r \in
\F_{q}[t_1,\ldots,t_n]$ be nonzero polynomials.   If $Z:=Z_{\F_{q}}(f_1,\dotsc,f_r)$,
then
\begin{equation*}
\ord_p(\# Z) \,\geq\, \bigg{\lceil} \frac{N\bigl(n-\sum_{j=1}^r\sigma_p(f_j)\bigr)}
{\max_{j=1}^r \sigma_p(f_j)} \bigg{\rceil}.
\end{equation*}
\end{thm}

\subsection{Functional degrees of polynomial functions in positive characteristic}
Let $R$ be a commutative ring of prime characteristic $p$.  Must we have equality in
\eqref{PWEIGHTINEQ}?  When $R$ is a field, this is answered by \cite[Thm.\
10.3]{Aichinger-Moosbauer21}.  In this result Aichinger-Moosbauer show that
$\fdeg(E(f)) = \sigma_p(f)$ whenever $R$ is an infinite field of characteristic $p$.
Later in this section we will show that this result continues to hold whenever $R$ is an
infinite \emph{domain} of characteristic $p$.
\\ \\
The case of $R = \F_q$ is more closely related to the main results in this paper:
a strict inequality
$\fdeg(E(f)) < \sigma_p(f)$ would yield a further improvement of the Ax-Katz
Theorem. It turns out that strict inequality can occur, however in a way that leads only
to improvements of the Ax-Katz Theorem that had already been well understood.
\\ \\
To explain, we call a nonzero monomial term $c_{\underline{d}}\, t_1^{d_1} \dotsm
t_n^{d_n} \in \F_q[t_1,\dotsc,t_n]$ \textbf{reduced} if $d_j \leq q-1$ for all $1 \leq j\leq
n$\,.  (Note the strong dependence on the ground field.)  A polynomial is
\textbf{reduced} if each of its nonzero monomial terms are reduced.  \\ \indent Just
using the fact that $x^q = x$ for all $x \in \F_q$, it is easy to see that to every $f \in
\F_q[t_1,\dotsc,t_n]$ there is a reduced polynomial $\overline{f} \in
\F_q[t_1,\dotsc,t_n]$ that induces the same function $\F_q^n \ra \F_q$ as $f$.
Already in \cite{Chevalley35}, Chevalley showed that every function in
$\F_q^{\F_q^n}$ is equal to $E(f)$ for a unique reduced polynomial $f$.  (For an
English language proof and some modest generalizations, see \cite[\S 2.3 and \S
3.1]{Clark14}.)  In particular, the polynomial $\overline{f}$ alluded to above is the
unique reduced polynomial inducing the same function as $f$.
\\ \\
Now for any $f_1,\dotsc,f_r \in \F_q[t_1,\dotsc,t_n]$, since the solution set
\[ Z(f_1,\dotsc,f_r) \,\coloneqq\, \{x \in \F_q^n \mid f_1(x) = \dotsb = f_r(x) = 0 \} \]
depends only on the associated functions $E(f_1),\dotsc,E(f_r)$, we always have
\[ Z(f_1,\dotsc,f_r) \,=\, Z(\overline{f_1},\dotsc,\overline{f_r}). \]
One gets easy strengthenings of many results of Chevalley-Warning type -- in
particular the theorems of Chevalley-Warning and Ax-Katz -- by replacing
$f_1,\dotsc,f_r$ by $\overline{f_1},\dotsc,\overline{f_r}$, since in this process none of
the degrees can increase.
\\ \\
The following result is part of \cite[Thm.\ 10.3]{Aichinger-Moosbauer21}.

\begin{thm}[Aichinger-Moosbauer]
\label{LASTAMTHM} Let $f \in \F_q[t_1,\dotsc,t_n]$ be a nonzero polynomial, and let
$E(f) \in \F_q^{\F_q^n}$ be the associated polynomial function. Then
\[ \fdeg(E(f)) \,=\, \sigma_p(\overline{f}) .\]
\end{thm}
\noindent Proposition \ref{PWEIGHTPROP} and Theorem \ref{LASTAMTHM} imply
that
\[ \sigma_p(\overline{f}) \,=\, \fdeg(E(f)) \,\leq\, \sigma_p(f)\,; \]
that is, passing to the reduced polynomial also cannot increase the $p$-weight
degree.  So, in the setting of the Moreno-Moreno Theorem, one can improve the
conclusion to
\begin{equation*}
\ord_p(\# Z_{\F_q^n}(f_1,\dotsc,f_r))
\,\geq\, \bigg{\lceil} \frac{N\bigl(n-\sum_{j=1}^r\sigma_p(\overline{f_j})\bigr)}
{\max_{j=1}^r \sigma_p(\overline{f_j})} \bigg{\rceil}
\end{equation*}
which by Theorem \ref{LASTAMTHM} is the optimal application of Corollary
\ref{GTPAKT} to polynomials over $\F_q$.

\begin{remark}
For a reduced polynomial $f \in \F_p[t_1,\dotsc,t_n]$, we have $\deg(f) =
\sigma_p(f)$, so Moreno-Moreno gives no essential improvement upon Ax-Katz when
$q = p$.
\end{remark}
\noindent Now we prepare for our generalization of 
\cite[Thm.\ 10.3]{Aichinger-Moosbauer21} with the following result.

\begin{lemma}
\label{FINALLEMMA} Let $A$ be an infinite commutative codomain that is a finitely
generated $\F_p$-algebra.  Let $x_1,\ldots,x_m \in A \setminus \{0\}$, and let $M \in
\Z^+$\!.  Then there is a maximal ideal $\mmm$ of $A$ such that:
\begin{itemize}
\item[(i)] $x_1,\ldots,x_m \notin \mmm$, and
\item[(ii)] $A/\mmm$ is a finite field of size greater than $M$.
\end{itemize}
\end{lemma}

\begin{proof}
Replacing $x_1,\ldots,x_m$ with $x \coloneqq x_1 \cdots x_m$, we reduce to the
case of $m = 1$.  Zariski's Lemma \cite[Thm.\ 11.1]{Clark-CA} implies that for every
maximal ideal $\mmm$ of $A$, the field $A/\mmm$ is a finite-dimensional
$\F_p$-vector space, hence a finite field.  The same argument shows that $A$ is not
itself a field; also $A$ is a Noetherian ring \cite[Cor.\ 8.39]{Clark-CA}.    Moreover,
since $A$ is a finitely generated algebra over the field $\F_p$ it is a Jacobson domain
\cite[Prop.\ 11.3b)]{Clark-CA}, so $\bigcap_{\mmm \in \MaxSpec A} \mmm = (0)$.   In
particular $A$ has infinitely many maximal ideals, since a finite intersection of
nonzero ideals in a domain is nonzero. It follows that the set $\mathcal{U}(x)$ of
maximal ideals $\mmm$ of $A$ such that $x \notin \mmm$ is nonempty.  We claim
that $\mathcal{U}(x)$ is moreover infinite: if on the contrary we had $\mathcal{U}(x) =
\{ \mmm_1,\ldots,\mmm_n\}$, then for $1 \leq i \leq n$ choose $y_i \in \mmm_i
\setminus \{0\}$, and we see that $xy_1 \cdots y_n$ is a nonzero element of $A$ that
lies in every maximal ideal of $A$: contradiction.  Finally, by \cite[Thm.\
22.23]{Clark-CA}, in any Noetherian ring $S$, for all $M \in \Z^+$ there are only
finitely man ideals $I$ of $M$ such that $S/I$ is finite of size at most $M$. So in any
infinite family of maximal ideals of $A$, the size of the residue ring approaches
infinity.
\end{proof}

\begin{thm}
\label{LASTTHM} Let $R$ be an infinite commutative domain of characteristic $p$.
Let $f \in R[t_1,\ldots,t_n]$ and let $E(f) \in \PP(R^n,R)$ be the associated polynomial
function. Then
\[ \fdeg(E(f)) = \sigma_p(f). \]
\end{thm}

\begin{proof}
By Proposition \ref{PWEIGHTPROP} it suffices to show that $\fdeg(E(f)) \geq \sigma_p(f)$.  Let $K$ be the fraction field of $R$. \\
\emph{Case 1, $K/\F_p$ is an algebraic field extension:}  This case is already
covered by \cite[Thm. 10.3]{Aichinger-Moosbauer21}, as necessarily $R = K$:
indeed, for every nonzero element $x \in R$ there is a positive integer $n_x$ such
that $x^{n_x} = 1$, so $x^{-1} = x^{n_x-1} \in R$.  \\
\emph{Case 2, $K/\F_p$ is transcendental:}  In this case $R$ must contain elements
that are transcendental over $\F_p$: let $t$ be such an element, and let $A$ be the
$\F_p$-subalgebra of $R$ generated by $t$ and the coefficients of $f$. Let $a \in
\Z^+$\!.  By Lemma \ref{FINALLEMMA}, there is a maximal ideal $\mmm$ of $A$ that
does not contain any of the coefficients of $f$ and such that $\F \coloneqq A/\mmm$
is a finite field $\F$ of order at least $p^a$.  Let $\underline{f}$ be the image of $f$ in
$\F[t_1,\ldots,t_n]$.  By Lemma \ref{HOMFUNCI} we have
\[ \fdeg(E(f)) \,\geq\, \fdeg(E(f)|_{A^n}) \,\geq\, \fdeg(E(\underline{f})). \]
By our choice of $\mmm$, the monomials appearing in $\underline{f}$ with nonzero
coefficient are the same as those appearing in $f$ with nonzero coefficient, so
$\sigma_p(\underline{f}) = \sigma_p(f)$.  Choosing $a$ larger than $\max_{1 \leq i
\leq n} \deg_i(f)$ makes $\underline{f}$ $\F$-reduced, so by Theorem
\ref{LASTAMTHM} we have
\[ \fdeg(E(f)) \,\geq\, \fdeg(E(\underline{f})) \,=\, \sigma_p(\underline{f})
\,=\, \sigma_p(f). \qedhere \]
\end{proof}

\section{Further Work}
\noindent It is natural to ask for a generalization of Theorem \ref{AXKATZWILSON} in
which instead of $(\Z/p\Z)^N$, we may take $A$ to be any finite commutative
$p$-group.  Such a result will be given in the forthcoming work \cite{Clark-Schauz23}.
The proof follows the same basic strategy: the $p$-adic divisibility comes from a
combination of Corollary \ref{COR4.57} and a generalization of Lemma \ref{Lem5.1}
to sums of the form $\int_{ \prod_{i=1}^N [p^{\alpha_i})} f$.  \\ \indent Here is a quick
overview of this work: to solve the number-theoretic problem of determining
$\ord_p(\int_{ \prod_{i=1}^N [p^{\alpha_i})} f)$ for
\[ f(\underline{x}) \,=\,\binom{x_1}{n_1} \dotsm \binom{x_N}{n_N} \]
in terms of $n_1,\ldots,n_N$ is not very difficult, but to solve the discrete optimization
problem of, for each fixed $\underline{\alpha} = (\alpha_1,\ldots,\alpha_N)$,
minimizing this quantity over all $\nn = (n_1,\ldots,n_N) \in \N^N$ with fixed $d
\coloneqq |\nn|$ takes more work.  Then we must minimize the total $p$-adic
divisibility obtained from this and from Corollary \ref{COR4.57}.  The answer obtained
is intricate in the general case, suggesting that these complications may be inherent
to the problem.
\\ \\
The $\beta = 1$ case of Lemma \ref{Lem5.1} gives a result of Ax \cite{Ax64}; let's call
it \emph{Ax's Lemma}.  Ax's Lemma comprises most of the ten line proof of
Chevalley-Warning referred to in the introduction.  It suggests a further problem in the
Aichinger-Moosberger calculus.

\begin{ques}
Let $A$ and $B$ be finite commutative $p$-groups.  What is the largest $d \in
\tilde{\N}$ such that for all $f \in B^A$ with $\fdeg(f) \leq d$, we have
\[ \int_A f \,\coloneqq\, \sum_{x \in A} f(x) \,=\, 0\,? \]
\end{ques}
\noindent Let us call this largest possible $d$ the \textbf{summation invariant}
$\sigma(A,B)$.  In this notation, Ax's Lemma amounts to:
\[ \forall p \in \mathcal{P}, \ \forall N \in \Z^+\!, \ \sigma\bigl( (\Z/p\Z)^N\!,\,\Z/p\Z\bigr)
 \,=\, N(p-1) - 1. \]
In the paper \cite{Clark-Triantafillou23} the following generalization of Ax's Lemma will
be shown:
\begin{equation}
\label{GENAXEQ}
 \forall p \in \mathcal{P}, \ \forall N \in \Z^+\!, \ \forall 1 \leq \beta \leq N, \
 \sigma( (\Z/p\Z)^N\!,\,\Z/p^\beta \Z) \,=\, N(p-1) -1.
\end{equation}
In \cite{Clark-Triantafillou23} we use (\ref{GENAXEQ}) to derive a qualitative
generalization of Ax-Katz over any finite rng $R$ of size divisible by $p$: if we fix the
number and degrees of polynomials $f_1,\ldots,f_r$, then $\ord_p(\#
Z_{R^N}(f_1,\ldots,f_r))$ approaches infinity with the number $N$ of variables.  \\
\indent Such results also follow from the main theorem of \cite{Clark-Schauz23} -- but
the argument is different.  Unlike the proofs presented here and in
\cite{Clark-Schauz23}, the arguments of \cite{Clark-Triantafillou23} do not use the
fundamental representation (Theorem \ref{4.2}): they work entirely in finite
characteristic.
\\ \\
Notice that (\ref{GENAXEQ}) is a finite characteristic variant of Lemma \ref{Lem5.1},
but when $\beta > 1$ the conclusion of (\ref{GENAXEQ}) is stronger than the
conclusion of Lemma \ref{Lem5.1}.  At first this seems strange: if in Lemma
\ref{Lem5.1} we weakened $\fdeg(f) < (p-1)(N-\beta + 1)$ to $\fdeg(f) < N(p-1)$ then
in general it is false that $\ord_p( \int_{[p)^N} f) \geq \beta$.  However, to compare the
two results we must take $f: (\Z/p\Z)^N\! \ra \Z/p^{\beta} \Z$ and pull it back to $F:
\Z^N\! \ra \Z/p^{\beta} \Z$, so the bound of (\ref{GENAXEQ}) applies only to functions
$F: \Z^N\! \ra \Z/p^{\beta} \Z$ that are $p$-periodic, whereas Lemma \ref{Lem5.1}
applies to all functions $F: \Z^N\! \ra \Z/p^{\beta} \Z$.  Thus in the application of
Lemma \ref{Lem5.1} to results on maps between finite commutative groups, we are
losing critical information, namely periodicity properties of the functions coming from
the fact that they were pulled back from finite characteristic.  This explains why our
present approach also includes Theorem \ref{THM4.57}, which uses the periodicity
properties to deduce further $p$-adic divisibilities coming from the coefficients of the
fundamental representation.
\\ \\
The two-pronged approach taken here and in \cite{Clark-Schauz23} seems to be
quantitatively superior to the approach via $\sigma(A,B)$ alone taken in
\cite{Clark-Triantafillou23}, but it would be interesting to clarify the relationship
between them.

\end{document}